\documentclass[11pt]{amsart}  
\usepackage{amscd,amssymb,amsmath,latexsym,enumerate}
\usepackage[mathscr]{euscript}
\usepackage{epsfig}

\newtheorem{theo}{Theorem}
\newtheorem{defini}{Definition}
\newtheorem{proposi}{Proposition}
\newtheorem{lemma}{Lemma}
\newtheorem{coro}{Corollary}
\newtheorem{rem}{Remark}
\newtheorem{exam}{Example}

\newtheorem{conj}{Conjecture}

\newcommand{\Ff}{{\mathcal F}}
\newcommand{\Gg}{{\mathcal G}}
\newcommand{\Hh}{{\mathcal H}}

\newcommand{\Ll}{{\mathcal L}}

\newcommand{\Pp}{{\mathcal P}}

\newcommand{\Uu}{{\mathcal U}}

\newcommand{\EM}{{\mathbb E}}

\newcommand{\NM}{{\mathbb N}}

\newcommand{\RM}{{\mathbb R}}

\newcommand{\ZM}{{\mathbb Z}}

\newcommand{\GG}{{\mathfrak G}}

\newcommand{\hFf}{{\widehat{\mathcal F}}}

\newcommand{\es}{{\mathscr E}}

\newcommand{\ts}{{\mathscr T}}
\newcommand{\us}{{\mathscr U}}
\newcommand{\vs}{{\mathscr V}}
\newcommand{\ws}{{\mathscr W}}

\newcommand{\diam}{\mbox{\rm diam}}                
\newcommand{\idiam}{\mbox{\rm \scriptsize diam}}   
\newcommand{\ch}{\mbox{\rm \tiny Ch}}	           
\newcommand{\pr}{\mbox{\rm Prob}}                  
\newcommand{\var}{\mbox{\rm Var}}                  

\begin{document}

\title{Bi-Lipshitz Embedding of Ultrametric Cantor Sets into Euclidean Spaces
\vspace{.5cm}}

\author{Jean V. Bellissard, Antoine Julien
\vspace{.2cm}
}

\address{Jean V. Bellissard, Georgia Institute of Technology, School of Mathematics, Atlanta GA 30332-0160
}

\email{jeanbel@math.gatech.edu}

\address{Antoine Julien, Department of Mathematical Sciences, NTNU, NO-7491, Trondheim, Norway
}

\email{antoine.julien@math.ntnu.no }
\thanks{Work supported in part by NSF grants DMS 0901514, DMS 1160962 and by the CRC701, {\em Spectral Structures and Topological Methods in Mathematics}, Mathematik, Universit\"at Bielefeld, Germany.
Substantial part of this work was made while the second author was employed by the University of Victoria, with the support of the Pacific Institute for the Mathematical Sciences.}

\begin{abstract} 
An ultrametric Cantor set can be seen as the boundary of a rooted weighted tree called the Michon tree. The notion of Assouad dimension is re-interpreted as seen on the Michon tree. The Assouad dimension of an ultrametric Cantor set is finite if and only if the space is bi-Lipschitz embeddable in a finite dimensional Euclidean space. This result, due to Assouad and refined by Luukkainen--Movahedi-Lankarani is re-proved in the Michon tree formalism. It is applied to answer the embedding question for some spaces which can be seen naturally as boundary of trees: linearly repetitive subshifts, Sturmian subshifts, and the boundary of Galton--Watson trees with random weights. Some of these give examples of nonembeddable spaces with finite Hausdorff dimension.
\end{abstract}


\maketitle

\section{Introduction}
\label{embed11.sect-intro}

\noindent This article is dedicated to giving a concrete description of the bi-Lipshitz embedding theorems for ultrametric Cantor sets in view of their applications to various situations like the case of tiling spaces or to the boundary of rooted weighted trees. In particular several important examples will be given of ultrametric Cantor sets with finite Hausdorff dimension that are not bi-Lipshitz embeddable in a finite dimensional Euclidean space.

\vspace{.1cm}

\noindent The problem of embedding metric spaces spaces in the standard Euclidean space $\RM^n$ is fairly old, and has many faces. Hence Fr\'echet's Theorem \cite{Fr10} shows that any finite metric space with $n+1$ elements can be embedded isometrically in $\RM^n$ when equipped with the $\ell^\infty$-norm. Schoenberg \cite{Sch37,Sch38} gave a necessary and sufficient condition for a finite metric space to be embeddable isometrically in the Euclidean space $\RM^n$. These works are commonly used in data analysis today. Beyond finite metric space, the case of compact ultrametric spaces is probably the simplest to consider next. An ultrametric space is a metric space where the distance satisfies the strengthened triangle inequality: $d(x,z) \leq \max \{ d(x,y),d(y,z) \}$ for any three points $x,y,z$ in the space. It is known \cite{ML92} that any such space is isometrically embeddable in a (real) Hilbert space.

\vspace{.1cm}

\noindent However the requirement that the embedding is isometric is too strong in general. For instance an ultrametric space is unlikely to be isometrically identified with a subset of a finite dimensional Euclidean space. The condition that the embedding should be isometric has to be weakened. A map $f: (X,d) \to (\RM^n, \| . \| )$ is a bi-Lipschitz embedding if it is a homeomorphism onto its image, and furthermore

$$\exists c\geq 1\,, 
   \hspace{1cm}
    \forall x,y\in X\,,
     \hspace{2cm}
     \frac{1}{c}\,d(x,y)\leq \|f(x)-f(y)\|\leq c\,d(x,y)\,.
$$

\noindent In the rest of this article, an \emph{embedding} will always be assumed to be bi-Lipschitz, and a space will be called embeddable (or $f$-embeddable, $f$ standing for ``finite-dimensional'') if it can be embedded into some $\RM^n$ equipped with its Euclidean metric.

\vspace{.1cm}

\noindent The concept of Assouad dimension has been a breakthrough in the problem of embeddability \cite{As77, As79, As83} (see also Semmes~\cite{Se99}). A space has finite Assouad dimension if and only if it is doubling, that is, for some $M\in\NM$, any ball can be covered by at most $M$ balls of half radius. Moreover a separable metric space is $f$-embeddable only if its Assouad dimension is finite. If, in addition, the space is ultrametric, then its $f$-embeddability is equivalent to having finite Assouad dimension. This result, obtained by Assouad, was strengthened by Luukkainen and Movahedi-Lankarani~\cite{LM94} and then Luosto~\cite{Lu96}: an ultrametric space is bi-Lipschitz embeddable in $\RM^n$ if and only if its Assouad dimension is strictly less than~$n$. The Assouad dimension seems therefore perfectly fitted to answer the embedding question for ultrametric spaces.

\vspace{.1cm}

\noindent In this paper, the focus is on ultrametric Cantor sets. Cantor sets are totally disconnected metrizable compact spaces, without isolated points. A theorem of Brouwer~\cite{Br10} states that any two such sets are homeomorphic (hence any of these sets can be called ``\emph{the} Cantor set''). However, the Cantor set can carry a lot of different ultrametrics. A classification of all ultrametric on a Cantor set was given by Michon~\cite{Mi85}: an ultrametric Cantor set can be represented by a rooted infinite tree, and the ultrametric is described by weights on the vertices (which correspond to the diameters of a basis of neighborhoods).

\vspace{.2cm}

\noindent This paper is structured as follows: Section~\ref{embed11-sect-defini} gives the basic notations and terminology. In particular the definition of Assouad dimension of a metric space and the description of ultrametric Cantor sets arising as boundaries of weighted rooted trees are described (which by Michon's Theorem, are all ultrametric Cantor sets in existence).

\vspace{.1cm}

Section~\ref{embed11.sect-embed-results} is dedicated to existing results about bi-Lipschitz embeddings. The interpretation of the Assouad dimension in terms of the Michon tree is discussed. Luukkainen--Movahedi-Lankarni's result is re-proved in this formalism and takes a very concrete form. The proof is relatively short, and self-sustained. It also has the advantage that the embedding seems slightly more explicit than in the original proof.

As a corollary (see \cite{JS11}), this result shows that the classical quasicrystalline tilings, like the Penrose tiling, the octagonal tiling and the three classes of icosahedral tilings used to model quasicrystalline metallic alloys \cite{LectQC} are embeddable when equipped with the combinatorial metric. Actually if $d$ is the dimension of such a tiling, the Hausdorff dimension of its tiling space is equal to $d$ and is embeddable in $\RM^{d+1}$, showing that the Assouad dimension is larger than or equal to $d$ and smaller than $d+1$. It is likely that, in such cases, the Assouad dimension is actually equal to $d$.

\vspace{.1cm}

Section~\ref{embed11.sect-sadic} gives a class of examples defined through a special type of subshift. An ultrametric Cantor set appears as the transversal of quasiperiodic tilings (or in one dimension, as subshifts of $\{0,1\}^{\ZM}$). The following two theorems are proved.

\begin{theo}
\label{embed11.th-LR}
The tiling space of a (bi-infinite, one-dimensional) linearly repetitive sequence  is $f$-embeddable. 
\end{theo}

Section~\ref{embed11.sect-sturm} investigates the case of Sturmian sequences leading to

\begin{theo}
\label{embed11.th-sturmHull}
The tiling space of a Sturmian sequence $x$ associated with the irrational number $\alpha\in (0,1)$ is $f$-embeddable if and only if $\alpha$ has bounded type, that is if and only if its continued fraction expansion has bounded partial quotients. In particular for Lebesgue-almost every $\alpha$ this tiling space is not $f$-embeddable.
\end{theo}

\noindent Interestingly enough the Minkowski (or box-counting) dimension of a subshift is related with its complexity (see \cite[Corollary 3.15]{JS11}). However, its $f$-embeddability (and hence its Assouad dimension) has more to do with the recurrence properties of the orbits, and therefore reflects some of the dynamical properties of the subshift. Both combinatorial and dynamical properties of the subshift seem to be captured by the metric. This remark is one of the motivations for studying the metric properties of ultrametric Cantor sets arising from tilings: the distance provides an additional structure which reflects dynamical or combinatorial properties of the tiling. A notion of spectral dimension associated with a family of spectral triples defined on this space had been investigated by Pearson and the first author~\cite{PB08}.
It can be seen here that each different concepts of dimension give a specific perspective on these spaces.

\vspace{.1cm}

In Section~\ref{embed11.sect-haus}, the Hausdorff dimension is discussed in terms of the Michon tree representation. In a previous paper, Savinien and the second author~\cite{JS11} showed that a \emph{self-similar} ultrametric Cantor set is embeddable in $\RM^n$ as soon as $n$ is strictly more than the Hausdorff dimension of the space. However, in general, the Hausdorff dimension has no reason to be a good predictor of the embeddability of a space. The Hausdorff dimension is also discussed to prepare the last example discussed in this paper.

\vspace{.1cm}

Section~\ref{embed11.sect-GW} concerns random trees obtained from a Galton--Watson branching process \cite{Ha63,AN72}. 
Let $p$ be a probability defined on the set $\NM=\{0,1,2,\cdots\}$ of natural integers. A $GW$-tree is defined as a rooted tree in which each vertex $v$ has $M_v$ children, where the $M_v$ are independent, identically distributed random variables with probability distribution $p$ ($p$ is a distribution on the natural integers).
It is a classical result in probability, partly proved initially by the Reverend Watson~\cite{GW1875}, called the {\em Galton--Watson--Haldane--Steffensen critical Theorem}~\cite{Ke66}, that if the average number of children $\langle M_v\rangle$ is less than or equal to one, the tree obtained in this way is almost surely finite (extinction).
In particular its reduced tree is empty. On the other hand, as was eventually proved by Steffensen in 1930, if $\langle M_v\rangle> 1$, the probability of extinction is less than one. However, the random tree produced in this way is likely to have dandling vertices. One way to avoid such a property is to force every vertex to have at least two children. This can be done by demanding that $p_0=p_1=0$, namely that $p$ be supported by $[2,\infty)\subset \NM$: then the Galton--Watson tree is automatically reduced and this model can be called a {\em reduced Random tree}. A first result is the following:

\begin{proposi}
\label{embed11.prop-GaltWatson}
Let $p$ be a probability on the set $[2,\infty)\subset \NM$. If it has an infinite support, the boundary of the reduced Galton--Watson tree associated with this probability is almost surely not $f$-embeddable.
\end{proposi}

\noindent In such a case a {\em random weight} can be added in the following way: a family $(\lambda_v)_{v\in\vs}$ of {\em i.i.d.} random variables supported in $[0,1]$, with common distribution $\rho$ is defined. A weight is defined inductively by setting that the weight of the root is equal to $1$, and $\kappa(v) = \lambda_v\kappa(u)$ is $v$ is a child of $u$. In order that this defines a weight, it is required that $\rho\{0\}=0$ and $\rho\{1\}<1$. Then 

\begin{theo}
\label{embed11.th-GWHaus}
Let $p$ be a probability on the set $[2,\infty)\subset \NM$ and let $\rho$ be a probability on $[0,1]$ such that $\rho\{0\}=0$ and $\rho\{1\}<1$. If $m=\langle M_v\rangle<\rho\{1\}^{-1}$, let $s=s_m$ be the unique solution of $\langle \lambda^s\rangle m=1$. Then the reduced random tree $\ts$ produced by the Galton--Watson process associated with $p$ and endowed with the random weight associated with $\rho$ gives rise to an ultrametric Cantor set $\partial\ts$ with Hausdorff dimension $s_m$ almost surely. In addition, its Hausdorff measure exists almost surely and is a random probability measure. 
\end{theo}

\noindent Note that the Hausdorff dimension had been computed previously when the random variables $\lambda_v$ are deterministic (for example equal to $e^{-1}$ or $1/2$), see Hawkes~\cite{Haw81}.

\vspace{0.3cm}

\noindent Several observations can be made as a result of the analysis carried out in this paper. First, the Sturmian and Galton--Watson trees provide very natural examples of spaces which have finite Hausdorff dimension, but infinite Assouad dimension. More specifically, given a weighted tree in which the weights decrease exponentially along the branches, the Hausdorff dimension of its boundary is related to its \emph{average} number of children per vertex, while its Assouad dimension is related to the \emph{maximum} number of children per vertex.

\vspace{.1cm}

\noindent As a final anecdote, real trees growing in nature satisfy the requirement that the length of the branches decreases exponential fast with the generation rank. Hence they behave like embeddable trees (in this case they are embedded in $\RM^3$) if the weight is a measure of the length of their branches. The previous results give strong constraints about the growth of real trees since the Assouad dimension is at most $3$. 

\vspace{.3cm}

\noindent {\bf Acknowledgements: }J.~B. thanks A. Grigor'yan and his group at the CRC701, University of Bielefeld for giving him the opportunity to present this result prior to publication. A.~J.\ thanks J.\ Savinien for discussions on this embedding problem.

\section{Definitions and notations}
\label{embed11-sect-defini}

 \subsection{Trees and ultrametric Cantor sets}
 \label{embed11-ssect-trees}

\noindent In this paper, a rooted tree $\ts$ is described by a set of vertices $\vs$, which is countable and satisfies the following properties:

(i) $\vs$ is partitioned into $\vs = \bigsqcup_{n \geq 0} \vs_n$, and a predecessor (or ``parent'') map is defined $\vs_{n+1} \rightarrow \vs_n$ for all $n \geq 1$,

(ii) $\vs_0$ consists of a unique element, called the root and noted $\bullet$. 

\noindent Given a vertex $v$, any vertex which admits $v$ as a parent is called a child of $v$. If $v \in \vs_n$, the integer $n$ is noted $|v|$ and called the \emph{depth} or the \emph{generation} of $v$.
There is a relation between vertices: $v \succeq w$ if there is a sequence of vertices $\gamma = (v_0, v_1, \ldots, v_n)$, with $v=v_0$ and $w=v_n$, such that each $v_k$ is the parent of $v_{k+1}$. This binary relation reads: ``$v$ is an ancestor of $w$'' or ``$w$ is a descendant of $v$''. It is an order relation with the root as maximal element. In this case, $\gamma$ is called a \emph{path} from $v$ to $w$, which \emph{goes through} vertices $v_0, v_1,\ldots, v_n$. A vertex is called {\em dangling} whenever it has no children, and {\em branching} if it has at least two children.

\vspace{.1cm}

\noindent It is sometimes convenient to encode the ``parent--child'' relation by means of edges (this is how trees are defined in graph theory): a tree can be given as a set of vertex and a set of edges $\es$, which are both partitioned as above, with two maps $r: \es_n \rightarrow \vs_{n+1}$ and $s:\es_n \rightarrow \vs_n$ (range and source respectively), such that $r$ is bijective. Then the parent map is simply $s \circ r^{-1}$. These definitions are equivalent (note in particular that a path can be defined by a sequence of edges or a sequence of vertices).

\vspace{.1cm}

\noindent The trees considered here are furthermore supposed to be locally finite (every vertex has finitely many children), and without dangling vertices. Under these assumptions, there is a well-behaved topological space, called the boundary of the tree, and defined as

$$\partial \ts = 
   \Big\{ (v_0, v_1, \ldots) \in \prod_{n \in \NM} \vs_n \,;\,
    \ v_0 = \bullet \,;\,
     \ v_k \text{ is the parent of } v_{k+1}  \Big\}
      \stackrel{def}{=}
      \lim_{\gets}(\vs_n,\to)\,,
$$

\noindent which is nothing but the definition of the inverse limit of the sets $\vs_n$ under the parent map $\vs_{k+1} \rightarrow \vs_{k}$. Under the local finiteness condition, all $\vs_n$ are finite. Therefore, equipped with the product topology, $\partial \ts$ is a totally disconnected compact set, as an inverse limit of finite sets. Given a vertex $v$, let $[v]$ be the set of all element $(v_n)_{n \in \NM}$ which pass through $v$. Then $\{ [v] \ ; \ v \in \vs \}$ is a basis for the topology. A rooted tree is called {\em Cantorian} whenever, in addition, each vertex admits a branching descendant. The boundary of a Cantorian tree has no isolated points, and is therefore a Cantor set. 

\vspace{.1cm}

\noindent A {\em weight} on the rooted tree $\ts$ is a map $\kappa:\vs\to (0,\infty)$ such that 

(i) if $w\preceq v$ then $\kappa(w) \leq \kappa (v)$,

(ii) $\lim_{|v|\to\infty} \kappa(v) =0$.

\noindent A weight $\kappa$ induces an ultrametric $d_\kappa$ on $\partial\ts$ by

$$d_\kappa(x,y) = \kappa(x\wedge y)\,,
   \hspace{2cm} x,y\in\partial\ts\,,
$$

\noindent where $x \wedge y$ will denotes the deepest vertex through which $x$ and $y$ pass.

\vspace{.1cm}

\noindent Conversely the Michon Theorem~\cite{Mi85} establishes that any ultrametric Cantor set is isometric to the boundary of a rooted Cantorian tree endowed with a weight. Such a tree is unique if it is additionally required to be \emph{reduced}, that is all of its vertices are branching. This tree is then called the Michon tree of the ultrametric Cantor set.

\vspace{.1cm}

\noindent Any given tree can be reduced as follows without changing its boundary (up to isometry).
Given $\ts$, let $\vs'\subset \vs$ be the set of branching vertices which have infinitely many branching descendants (if $\ts$ is Cantorian, $\vs'$ is the set of all branching vertices).
If $w \in \vs'$, define its parent $v$ in $\ts'$ to be the closest ancestor in $\ts$ which belongs to $\vs'$. The new root $\bullet'$ is either $\bullet$ if the latter belongs to $\vs'$ or its closest descendant belonging to $\vs'$. The vertex set $\vs'$ with the ancestor relation inherited from $\ts$ defines a tree $\ts'$, which is reduced.
If the original tree comes with a weight $\kappa$, then the restriction of $\kappa$ to $\vs'$ gives a weight on the reduced tree. The two boundaries endowed with the induced ultrametrics are then isometrically homeomorphic.

\vspace{.1cm}

\noindent In the rest of this paper, each time an ultrametric Cantor set $(X,d)$ is given, it is implied that the unique Michon tree which represents it $(\ts,\kappa)$ is given as well.

 \subsection{Assouad dimension}
 \label{embed11.ssect-assdim}

\noindent Assouad defined \cite{As77,As79} a metric dimension associated with a metric space. When the space is ultrametric, this dimension captures particularly well the extent to which the space can be embedded in an Euclidean space.

\begin{defini}
\label{embed11.def-assdim}
Given $s>0$, a metric space $(X,d)$ is $s$-subhomogeneous if there exist $C > 0$ such that, for all $a > 0$ and all $b > 0$, and for all discrete subset $Y \subset X$, 

$$\Big( \forall x \neq y \in Y, \ a \leq d(x,y) \leq b \Big) 
   \hspace{1cm}\Rightarrow\hspace{1cm}
     \# Y\leq C \left( \frac{b}{a} \right)^s.
$$

\noindent Then, the metric dimension (or Assouad dimension) of $X$ is defined by

$$\dim_\mathrm{A} (X)=
   \inf\{s\in [0,+\infty]\,;\, X\text{ is $s$-subhomogeneous}\}
$$
\end{defini}

\noindent It is important to notice that the definition is uniform over the space $X$. Namely both constants $C$ and $s$ are independent on the location of the subset $Y$. It becomes clear that $\dim_\mathrm{A} (X)$ is finite if and only if $(X,d)$ has the doubling property.

\begin{proposi}
\label{embed11.prop-biLipinv}
The Assouad dimension is invariant under bi-Lipschitz homeomorphism.
\end{proposi}

\section{Embedding Theorems}
\label{embed11.sect-embed-results}

 \subsection{Known Results}
 \label{embed11.ssect-known}

\noindent For ultrametric Cantor sets the following Assouad Theorem holds

\begin{theo}[Assouad~\cite{As77, As79}]
\label{embed11.th-assouad}
Let $(X,d)$ be an ultrametric space. If $\dim_\mathrm{A} (X,d) < + \infty$, then $(X,d)$, is bi-Lipschitz embeddable in an Euclidean space.
\end{theo}

\noindent This result was later refined by Luukkainen and Movahedi-Lankarani \cite{LM94}

\begin{theo}
\label{embed11.th-embedding}
Let $(X,d)$ be an ultrametric space. If $\dim_\mathrm{A} (X) < n \in \NM$, then $X$ can be embedded in $\RM^n$. Conversely, if $X$ can be embedded in $\RM^n$, then $\dim_\mathrm{A} (X) \leq n$.
\end{theo}

\noindent This last theorem provides an optimal dimension for embedding when the Assouad dimension of the space is not an integer. However, if $\dim_\mathrm{A} (X)$ is an integer, the theorem does not answer the question whether or not $X$ is embeddable in $\RM^n$. This is addressed by Luosto \cite{Lu96}, and the answer is negative

\begin{theo}[Corollary 4.6 in~\cite{Lu96}]
\label{embed11.th-Lu96}
An ultrametric space $(X,d)$ can be bi-Lipschitz embedded in $\RM^n$ if and only if $\dim_\mathrm{A} (X) < n$.
\end{theo}

 \subsection{Assouad Dimension and Michon Trees}
 \label{embed11.ssect-AdimTree}

\noindent The purpose of this Section is to give a description of the Assouad dimension for ultrametric Cantor sets in terms of its Michon tree representation.

\vspace{.1cm}

\noindent Given a weighted Cantorian tree $(\ts, \kappa)$, a vertex $v$ of the tree and a number $0 < \delta < 1$, the \emph{sub-tree of $\ts$ under $v$ at resolution $\delta$} (noted $\ts(v,\delta)$) is defined as follows.

\begin{itemize}
 \item Its vertices consist of all $w \in \vs$ such that $w \preceq v$ and $\kappa(w) \geq \delta \kappa (v)$, as well as their children;
 \item The root is $v$;
 \item Its edges are all edges of $\ts$ between vertices of the sub-tree.
\end{itemize}

\noindent Since $\ts$ is locally finite and the weight tends to $0$ at infinity, this sub-tree is finite. By definition, a \emph{leaf} (or minimal vertices for $\succ$) of this finite tree is a vertex $w$ in the subtree satisfying $\kappa (w) < \delta \kappa(v)$. Hence leaves are the \emph{closest descendants of $v$} satisfying this relation. Let $\Ll (\ts,v,\delta)$ denote the set of all such leaves.

\begin{proposi}
\label{embed11.prop-shomMichon}
Let $(X,d)$ be an ultrametric Cantor set with Michon tree $(\ts,\kappa)$.
Then $X$ is $s$-subhomogeneous if and only if there exists a constant $C > 0$ such that for all vertex $v$ in $\ts$ and all $0 < \delta < 1$,
\begin{equation}
\label{embed11.eq-Linq}
\# \,\Ll (\ts,v,\delta) \leq
   C\, \delta^{-s}\,.
\end{equation}
\end{proposi}

\begin{rem}
\label{embed11.rem-Munb}
{\em Since the number of leaves in $\ts(v,\delta)$ is bounded from below by the number of children of $v$, it follows immediately that the Assouad dimension is infinite whenever the number of children per vertex is not uniformly bounded.
}
\hfill $\Box$
\end{rem}
\begin{proof}(i) Let $(X,d)$ be $s$-subhomogeneous. For any vertex $v$ and any $0 < \delta < 1$, let $Y$ be the set made of exactly one point in each $[w]$ when $w$ runs through the set of leaves of $\ts(v,\delta)$. Since $Y \subset [v]$ it follows that $d(x,y) \leq \kappa (v)$ for $x,y \in Y$. Moreover, if $x\neq y$, then $d(x,y) = \kappa (u)$ for some vertex $u\in\ts(v,\delta)$ which is an ancestor of the leaves in which $x$ and $y $ are chosen. Hence $d(x,y) \geq \delta \kappa (v)$. The $s$-homogeneity of $(X,d)$ then gives the inequality~(\ref{embed11.eq-Linq}).

\vspace{.1cm}

\noindent (ii) Conversely, let $(X,d)$ satisfies~(\ref{embed11.eq-Linq}). Let $0 < a < b$, and let $Y$ be a subset of $X$ such that $a \leq d(x,y) \leq b$ for all $x,y \in Y$. In particular the diameter of $Y$ is less than or equal to $b$. Hence there is a vertex $v$ such that $Y\subset [v]$ with $\kappa(v)\leq b$. With each $x\in Y$ is associated the leaf $w_x$ of $\ts(v,a/b)$ through which the path representing $x$ is going. If $y\neq x$ is another point of $Y$, it follows that $w_y\neq w_y$. For otherwise $d(x,y)\leq \kappa(w_x)<\kappa(v)\,a/b\leq a$. Hence the cardinality of $Y$ cannot be larger than the cardinality of $\Ll(\ts,v,a/b)$. Thanks to eq.~(\ref{embed11.eq-Linq}), it follows that $(X,d)$ is $s$-homogeneous.
\end{proof}

\vspace{.2cm}

\noindent A criterion for checking finiteness of the Assouad dimension is given now. Remark that in conjunction with Theorem~\ref{embed11.th-assouad}, this gives a criterion for embeddability into some $\RM^n$.

\begin{proposi}
\label{embed11.prop-criterion}
Let $(X,d)$ be an ultrametric Cantor set with Michon tree $(\ts,\kappa)$.
 Then the Assouad dimension of $X$ is finite if and only if the number of children per vertex is bounded (uniformly) and there are two constants $c \geq 1$ and $0 < \alpha < 1$ such that for all vertices $v \succeq w$,

\begin{equation}
\label{embed11.eq-expdec}
\frac{\kappa (w)}{\kappa (v)} \leq c \, \alpha^{|w| - |v|}\,,
\end{equation}

\noindent  where $|v|$ (resp.\ $|w|$) is the depth of $v$ (resp.\ $w$) in the tree.
\end{proposi}

\begin{rem}
\label{embed11.rem-doub}
{\em This characterization in terms of the Michon tree can be compared with the doubling property for a metric space. In the present case, it corresponds to the condition that the cardinality of the sets $\Ll(\ts,v,1/2)$ is uniformly bounded. Assouad proved that the doubling property is equivalent to finite metric dimension.
}
\hfill $\Box$
\end{rem}
\begin{proof}(i) Let~(\ref{embed11.eq-expdec}) hold. Then let $M$ be the maximal number of children per vertex. For any $\delta > 0$, let $k-1$ be the integer part of $\ln (\delta) / \ln (\alpha)$. If $w$ is a descendant of $v$ of depth $|v|+k$, then

\begin{align*}
 \kappa (w)  &  \leq  c \, \kappa (v) \, \alpha^{k} 
     \leq \kappa (v) \exp\Big( \ln(\alpha) \frac{\ln(\delta)}{\ln(\alpha)} \Big)  = \kappa (v) \, \delta.
\end{align*}

\noindent Therefore, elements in $\Ll(\ts,v,\delta)$ are at depth at most $|v| + k$ and its cardinality is at most $M^k$. Consequently,

\begin{align*}
  M^k & \leq M \cdot M^{\ln(\delta) / \ln (\alpha)} 
       = M \exp{\Big( \ln (\delta) \frac{\ln(M)}{\ln (\alpha)} \Big)}  = M \delta^s,
 \end{align*}

\noindent with $s=\ln(M) / \ln (\alpha)$. Therefore, $X$ is $s$-subhomogeneous, and the Assouad dimension is finite.

\vspace{.1cm}

\noindent (ii) If the assumption represented by eq.~(\ref{embed11.eq-expdec}) does not hold then
\(
\forall C \geq 1, \  \forall 0 < \alpha < 1, \ \exists w \preceq v
\)
vertices in $\ts$, such that
\[
 \frac{\kappa (w)}{\kappa (v)} > C \alpha^{|w| - |v|}.
\]

\noindent In particular it implies that

\begin{equation}
\label{embed11.eq-ndp}
\forall n \in \NM, \ 
 \exists v_n \succ w_n  
  \text{ such that } |w_n|-|v_n| = n 
   \text{ and }
    \frac{\kappa (w_n)}{\kappa (v)} > \frac{1}{2}.
\end{equation}

\noindent For indeed, if not, let $n$ be chosen such that for any vertex $v$ and any descendant $w$ of $v$ which is $n$ levels deeper, $\kappa (w) / \kappa (v) \leq 1/2$. Let now $v$ and $w$ be two arbitrary vertices such that $w \prec v$. By Euclidean division, $|w|-|v| = kn + r$ with $0 \leq r < n$. Then, let $v = w_0, w_1, \ldots, w_k$ be vertices between $v$ and $w$ such that $|w_i| - |v| = in$. By assumption,

\[
 \frac{\kappa (w)}{\kappa(v)} \leq \frac{\kappa (w_k)}{\kappa (v)}
    \leq \prod_{i=1}^{k} \frac{\kappa (w_i)}{\kappa(w_{i-1})} \leq \frac{1}{2^k}
    \leq 2 \left( \frac{1}{2^{1/n}} \right)^{kn+r} = \left( \frac{1}{2^{1/n}} \right)^{|w|-|v|}.
\]

\noindent So using the negation of equation~\eqref{embed11.eq-ndp}, it was possible to establish a geometric rate of decay of the weights, which was assumed not to hold in the first place. Therefore, equation~\eqref{embed11.eq-ndp} holds.

\vspace{.1cm}

\noindent (iii) If the assumption represented by eq.~(\ref{embed11.eq-expdec}) does not hold then, using equation~\eqref{embed11.eq-ndp}, there exists a sequence of vertices $v_n$ such that $\Ll(\ts, v_n, 1/2)$ contains $w_n$, which is $n$ levels deeper than $v_n$. So, the subtree $\ts (v_n, 1/2)$ contains all vertices between $v_n$ and $w_n$, as well as all their children. Since each vertex is branching, it has at least two children. Therefore, $\Ll(\ts, v_n, 1/2)$ contains at least $n$ elements.
This quantity tends to infinity as $n$ does. Hence, using the characterization of $s$-subhomogeneity of Proposition~\ref{embed11.prop-shomMichon}, $X$ cannot be $s$-subhomogeneous for any $s$, and its Assouad dimension is infinite.
\end{proof}

 \subsection{Embedding Theorems and Michon trees}
 \label{embed11.ssect-pfemb}

\noindent This section provides a proof of Theorem~\ref{embed11.th-embedding} for Cantor sets, using the Michon tree formalism.

\begin{proposi}
\label{embed11.prop-embedNec}
Let $(X,d)$ be an ultrametric Cantor set with Michon tree $(\ts,\kappa)$. If $(X,d)$ is bi-Lipschitz embeddable in $\RM^n$, then the Assouad dimension of $(X,d)$ is at most $n$.
\end{proposi}
\begin{proof}
\noindent Suppose $h$ is a bi-Lipschitz map $(X,d) \longrightarrow \RM^n$. Then, there are two constants $0 < m < M$ such that for any vertex $v$ of $\ts$,
\[
 m \kappa (v) \leq \diam ( h([v] ) \leq M \kappa (v).
\]

\noindent The goal is to show that $(X,d)$ is $n$-subhomogeneous, so that the Assouad dimension of $(X,d)$ is at most equal to $n$.
Let $0 < \delta < 1$, and let $v$ be a vertex of $\ts$. By definition of $\Ll(\ts,v,\delta)$, given any distinct $w,w' \in \Ll(\ts,v,\delta)$ and any $x \in w$ and $x' \in w'$, one has $d(x,x') = \kappa (w \wedge w') \geq \delta \kappa (v)$.
Therefore,
\[
 d(h(x),h(x')) \geq m \, d(x,x') \geq m \delta \kappa (v).
\]
\noindent Therefore, the Euclidean balls of center $h(x)$ and $h(x')$ and diameter $m \delta \kappa (v)$ are disjoint. It holds for any pair $w \neq w'$ of elements in $\Ll(\ts, v, \delta)$. It means that $h([v])$, which is included in a ball of diameter $M \kappa (v)$ contains $\# \Ll (\ts,v,\delta)$ disjoint balls of diameter $m \delta \kappa (v)$.
By a comparing the volumes, one gets:
\[
 \# \Ll (\ts,v,\delta) \big(m\delta\kappa(v)\big)^n \leq \big( M \kappa(v) \big)^n.
\]
\noindent Rearranging the terms, it gives $n$-subhomogeneity with constant $C=(M/n)^n$.
\end{proof}

\begin{proposi}
\label{embed11.prop-embedSuff}
Let $(X,d)$ be a $s$-subhomogeneous space, with $s < n$. Then $(X,d)$ is bi-Lipschitz embeddable in $\RM^n$.
\end{proposi}
\noindent The proof is adapted from the one in~\cite{JS11}. It relies on the construction of the $\delta$-compressed tree ($0 < \delta < 1$) in an way analogous to the construction of telescoped diagrams in the self-similar case. The construction is given by the following lemma.

\begin{lemma}
\label{embed11.lem-telesc}
Let $(\ts,\kappa)$ be a weighted Cantorian reduced tree, which is $s$-subhomogeneous. Then, for all $0<\delta<1$, there exists $(\ts_\delta, \kappa_\delta)$ a weighted Cantorian reduced tree called the $\delta$-compression of $(\ts,\kappa)$ such that

 \begin{enumerate}
  \item $\partial \ts$ and $\partial \ts_\delta$ are bi-Lipschitz homeomorphic (for the distances induced by their respective weight functions);
  \item $\kappa_\delta(w) < \delta \kappa_\delta (v)$ for all parent--child pair $(v,w)$ of vertices of $\ts_\delta$.
 \end{enumerate}

\noindent Furthermore, there exists a $C$ (independent of $\delta$), such that for all $\delta$, the number of children per vertex in $\ts_\delta$ is bounded uniformly by $C \delta^{-s}$.
\end{lemma}
\begin{proof}
\noindent (i) The $\delta$-compression of $\ts$ is built inductively as follows: let $\delta$ be chosen in $(0,1)$.

\begin{itemize}
 \item[(i)] The root of $\ts_\delta$ is $\bullet$, so that $(\vs_\delta)_0 = \vs_0$.

 \item[(ii)] $(\vs_\delta)_{n+1} = \bigcup\{ \Ll(\ts,v,\delta) \ ; \ v \in (\vs_\delta)_n \}$.

 \item[(iii)] $w \ in (\vs_\delta)_{n+1}$ is a child of $v \in (\vs_\delta)_{n}$ in $\ts_\delta$ in $\ts_\delta$ if and only if $w$ is a descendant of $v$ in $\ts$.

 \item[(iv)] The weight $\kappa'$ on $\vs'$ is the restriction of $\kappa$.
\end{itemize}

\noindent The new tree $\ts_\delta$ is always reduced if $\ts$ is, because a vertex $v$ has more children in $\ts_\delta$ than it has in $\ts$, and therefore all vertices of $\ts_\delta$ are branching. Also, by definition of $\Ll(\ts,v,\delta)$, the weight satisfy condition (2) of the lemma.

\noindent Thanks to the $s$-subhomogeneity, there is $C>0$ such that the number of children of a vertex in $\ts_\delta$ is bounded by the cardinality of the sets $\Ll(\ts,v,\delta)$ namely by $C \delta^{-s}$. It proves the last point. Another consequence is that $\ts_\delta$ is locally finite, and therefore Cantorian.

\vspace{.1cm}

\noindent (ii) The map $\phi:\partial \ts\to\partial \ts_\delta$ is defined as follows. Let $\gamma \in \partial \ts$ be a path defined by the  sequence of vertices $(v_0, v_1, \ldots)$ be the sequence of vertices through which it passes. Since $\lim_{k\to]infty}\kappa (v_k)=0$, there is a sequence $(v_{k(i)})_{i \geq 0}$ such that $k(0)=0$ and $k(i+1)=\min\{k>k(i)\,;\, \kappa (v_{k(i+1)}) < \delta \kappa(v_{k(i)})\}$. Hence the $v_{k(i)}$ are vertices of $\ts_\delta$. Consequently the sequence of vertices $(v_{k(i)})_{i \geq 0}$ defines a unique path in $\ts_\delta$, which is set to be $\phi (\gamma)$ by definition. It is easy to see that $\phi$ is bijective, and that the image by $\phi$ of an open set is open. Therefore, $\phi^{-1}$ is continuous. Since $\ts_\delta$ is locally finite (by the $s$-subhomogeneity assumption), $\partial \ts_\delta$ is compact, so $\phi$ is a homeomorphism.

\vspace{.1cm}

\noindent (iii) To see that $\phi$ is bi-Lipschitz, let $x,y \in \partial \ts$. Let $v = x \wedge y$. Then $d(x,y) = \kappa (v)$. Let $w$ be the closest ancestor of $v$ which is also a vertex of $\ts_\delta$ (possibly $v=w$; more precisely, $v$ is a vertex of $\ts_\delta$ if and only if $v=w$). Then by definition of $\ts_\delta$,

\begin{equation}
\label{embed11.eq-bilip}
\delta \kappa (w) \leq \kappa (v) \leq \kappa (w).
\end{equation}

\noindent It is a quick check that $d_{\kappa_\delta} (\phi(x),\phi(y)) = \kappa_\delta (w)$. Together with equation~\eqref{embed11.eq-bilip}, it proves that $\phi$ is bi-Lipschitz.
\end{proof}

\begin{proof}[Proof of Proposition \ref{embed11.prop-embedSuff}]
\noindent Let $(\ts,\kappa)$ be the Michon tree associated with $(X,d)$. Let $0 < \delta < 1$, and $(\ts_\delta,\kappa_\delta)$ be the $\delta$-compression of the tree. The value of $\delta$ will be adjusted later.
The goal is to first build a map $\partial \ts_\delta \rightarrow \RM^n$, and then show that this map is a bi-Lipschitz homeomorphism for an appropriate value of $\delta$. Since $(X,d)$ is bi-Lipschitz homeomorphic to $\partial \ts_\delta$, it will prove the result.

\noindent (i) Let $M(\delta)$ be an upper bound for the number of children per vertex in $\ts_\delta$ (there is a $C$ such that for all $\delta$, $M(\delta) \leq C \delta^{-s}$, with $s < n$).
Then, define a function $g: \vs_\delta \rightarrow \{1, \ldots, M(\delta)\}$ such that its restriction to the set of children of $v$ be one-to-one for all vertex $v\in\vs_\delta$. Such a map is simply a numbering of children of each vertex. Any point $x\in\partial\ts_\delta$ will be represented by the sequence $x=(v_k)_{k \geq 0}$ of its vertices in $\ts_\delta$, with  $v_k$ a vertex of depth $k$, for all $k$. In particular (by definition of $\ts_\delta$), $\kappa_\delta(v_k) < \delta^k$. Then let $\phi: \partial \ts_\delta \rightarrow \RM^n$ be defined by $\phi(x) = (\phi_r(x))_{r=1}^n$ with

$$\phi_r(x)=
   \sum_{j=0}^\infty g(v_{jn+r}) \kappa_\delta (v_{jn+r-1}) ,
\qquad
 1 \leq r \leq n.
$$

\noindent Since $g(v_k)\leq M(\delta)$ and $\kappa_\delta (v_k) < \delta^k$ it follows that the series converges absolutely and uniformly with respect to $r$ and $x$.

\vspace{.1cm}

\noindent (ii) Let now $x,y\in\partial\ts$ with $x\neq y$ with $x=(v_k)_{k \geq 0}$ and $y=(w_k)_{k \geq 0}$. Then let $j_0,r_0$ be natural integers such that $|x \wedge y| = j_0 n + r_0 - 1$ and $1 \leq r_0 \leq n$.
It follows that the paths $x$ and $y$ share the same vertices until the generation $|x \wedge y|$ and they split at the next generation. Therefore (with the notation that $\chi(r\geq r_0)$ is $1$ if $r \geq r_0$ and $0$ otherwise),
\begin{multline*}
\phi_r(x)-\phi_r(y) =
   \sum_{j=j_0 + 1}^\infty \bigl(
     g(v_{jn+r})\kappa(v_{jn+r-1})-g(w_{jn+r})\kappa(w_{jn+r-1})
                         \bigr)\\
+ \chi(r\geq r_0)\bigl(
     g(v_{j_0n+r})\kappa(v_{j_0n+r-1})-g(w_{j_0n+r})\kappa(w_{j_0n+r-1})  \bigr)
\end{multline*}

\noindent It should be remarked that 
\[
-g(w_{jn+r})\kappa(w_{jn+r-1})\leq
   g(v_{jn+r})\kappa(v_{jn+r-1})-g(w_{jn+r})\kappa(w_{jn+r-1})\leq 
    g(v_{jn+r})\kappa(v_{jn+r-1}).
\]
Furthermore, both $v_{jn+r-1}$ and $w_{jn+r-1}$ are descendants of $x\wedge y$, sitting $n(j-j_0) + (r-r_0)$ levels deeper than $x \wedge y$.
Therefore, by construction of $\ts_\delta$,
\[
 \kappa(v_{jn+r-1}) < \delta^{n(j-j_0) + (r-r_0)} \kappa (x \wedge y), \qquad \text{and similarly for } w_{jn+r-1}.
\]
\noindent It leads to the inequality
\[
|g(v_{jn+r})\kappa(v_{jn+r-1})-g(w_{jn+r})\kappa(w_{jn+r-1})|
   \leq
    M(\delta) \kappa(x\wedge y) \delta^{n(j-j_0) + (r-r_0)}
\]
\noindent Therefore
\begin{equation}
\label{embed11.eq-deltaphi}
 |\phi_r(x)-\phi_r(y)|\leq 
   \frac{M(\delta)}{1-\delta^n}\;\kappa(x\wedge y)\,.
\end{equation}
\noindent Since $\kappa(x\wedge y)=d_\kappa(x,y)$, this map is Lipschitz continuous.
On the other hand, $\|\phi(x)-\phi(y)\| \geq |\phi_r(x)-\phi_r(y)|$ for all $1\leq r \leq n$. In particular, using again eq.~(\ref{embed11.eq-deltaphi}) and since $v_{j_0 n + r_0} \neq w_{j_0 n + r_0}$, it follows that
\[
\|\phi(x)-\phi(y)\| \geq 
   |\phi_{r_0}(x)-\phi_{r_0}(y)|\geq 
    \kappa(x\wedge y) -
     \frac{\delta^n M(\delta)}{1-\delta^n}\;\kappa(x\wedge y)\,.
\]
\noindent Now, since $M(\delta)$ grows at most like $C \delta^{-s}$, with $s < n$, it follows that for a small enough value of $\delta$, this map is actually bi-Lipschitz (in particular, it is one-to-one).
\end{proof}

\section{$S$-adic Systems}
\label{embed11.sect-sadic}

 \subsection{Definitions}
 \label{embed11.ssect-def}

\noindent This section is devoted to the basic definitions for $S$-adic systems. The presentation uses the notations of Durand~\cite{Dur1,Dur2}. Let $A$ be a finite alphabet. Note $A^\star$ the set of all finite words with letters
in $A$. For $w\in A^\star$, $|w|$ denotes its length, namely the number of letters it contains. Consider $S$ a \emph{finite} set of morphisms

$$\sigma \in S : A(\sigma) \longrightarrow A^\star, \quad 
   \text{with } A(\sigma) \subset A.
$$

\noindent A $S$-adic system is a sequence $( \sigma_n : A_{n+1} \rightarrow (A_n)^\star)_{n \in \NM} \in S^{\NM}$, such that the morphisms are composable. It will be assumed that for all $n$, every letter in $A_n$ appears in a word $\sigma_n(b)$ for some $b \in A_{n+1}$. For $m > n$ the following notation will be used

$$\sigma_{n,m} = 
   \sigma_n \circ \ldots \circ \sigma_{m-1}: A_m \rightarrow A_n.
$$

\noindent It satisfies $\sigma_{n,m} \circ \sigma_{m,k} = \sigma_{n,k}$.
A $S$-adic system is \emph{primitive} if there exists some $s_0>0$ such that for all $r \in \NM$, for all $a \in A_{r+s_0}$ and all $b \in A_r$, the letter $b$ appears in the word $\sigma_{r,r+s_0} (a)$.
It is called \emph{proper} if there exist two letters $l$ and $r$ in $A$ such that for all  $\sigma \in S$ and all $a \in A(\sigma)$, $\sigma(a)$ begins by the letter $l$ and ends by $r$. Furthermore, it will be assumed from this point on, that

\begin{equation}
\label{embed11.eq-length}
 \lim_{n \rightarrow +\infty}
  \min_{c \in A_n} |\sigma_{1,n} (c)| = +\infty,
\end{equation}
so that words of arbitrary long length are obtained from iterating the substitutions on a single letter.

\noindent Given a proper $S$-adic system, there is a way to associate a subshift $\Xi \subset A^\ZM$ (it is then called a $S$-adic subshift). Let $T$ be the shift operator on $A^\ZM$

$$T(\ldots x_{-1} \cdot x_0 x_1 \ldots) = 
   \ldots x_{-1} x_0 \cdot x_1 \ldots.
$$

\noindent It is straightforward to extend the morphisms of $S$ to $A^\ZM$ by concatenation

$$\sigma (\ldots x_{-1} \cdot x_0 x_1 \ldots) = 
   \ldots \sigma(x_{-1}) \cdot \sigma(x_0) \sigma(x_1) \ldots.
$$

\noindent It ought to be remarked that, by properness, for all $n$, $\sigma_{1,n}(l)$ is a prefix of $\sigma_{1,n+1}(l)$. Similarly, $\sigma_{1,n}(r)$ is a suffix of $\sigma_{1,n+1}(r)$. Therefore, an element of $A^\ZM$ can be defined by

$$x = \bigg( \lim_{n \rightarrow +\infty} \sigma_{1,n}(r)\bigg) 
\cdot \bigg( \lim_{n \rightarrow +\infty} \sigma_{1,n}(l)
      \bigg),
$$

\noindent where the dot separates the $x_i$ with $i \geq 0$ on the right from the ones with $i<0$ on the left. The subshift associated with the $S$-adic system is the closure in $A^\ZM$ (endowed with the product topology) of the orbit of $x$ under the shift. This subshift, $\Xi$, is endowed with the combinatorial distance~$d$, namely two sequences have distance $(n+1)^{-1}$ whenever they coincide on a string of radius $n$ around the origin and do not coincide beyond. If the $S$-adic system is primitive, $(\Xi,T)$ is minimal (see Durand~\cite[Lemma~7]{Dur1}).

\begin{defini}
\label{embed11.def-LR}
A word $x \in A^\ZM$ is called linearly recurrent or linearly repetitive  (LR for short) if there is a constant $K$ such that for all $n \in \NM$, for all subwords $w$ and $w'$ of $x$ respectively of length $n$ and $Kn$, then $w$ occurs in $w'$ as a subword.
\end{defini}

\noindent In other terms, in a LR word, each finite subword repeats infinitely often and within a distance which varies linearly with its size. It is easy to see that a subshift generated by an LR word is minimal and that every elements of the subshift is LR with the same constant. In this case, the subshift itself is called {\em linearly recurrent}. The following is due to Durand.

\begin{theo}[Durand~\cite{Dur2}, Proposition~1.1]
\label{embed11.th-LRdur}
A subshift is $S$-adic primitive and proper if and only if it is linearly recurrent.
\end{theo}

\noindent The periodic case is trivial. In the non-periodic case, Durand's construction of the $S$-adic system associated with a linearly recurrent
subshift involves return words. The $S$-adic system he builds in this case has the unique decomposition property (defined below).
This property will be used later on (see~\cite[Section~4]{Dur2} for the construction, and~\cite[Definition~9, Lemma~17]{DHS} for the properties of return words and codes).

\begin{defini}
\label{embed11.def-udp}
Let $(\sigma_n)_{n \in \NM}$ be a $S$-adic system. It is said to have
unique decomposition property if, given any element $x$ in the subshift associated with $\Xi$, there is a unique decomposition of $x$ as a concatenation

$$x = \ldots \sigma_1 (a_{-1}) \underline{\sigma_1 (a_0)} \sigma_1 (a_1) \ldots\,,
$$

\noindent where the index $0$ of $x$ is in the underlined word $\sigma_1 (a_0)$.
\end{defini}

 \subsection{Proof of Theorem~\ref{embed11.th-LR}}
 \label{embed11.ssect-pfLR}

\noindent The periodic case is trivial (in this case, $\Xi$ is finite), and so it may be assumed that the subshift is aperiodic.
Using the theorem of Durand cited above, it is enough to prove the following

\begin{proposi}
\label{embed11.prop-primEmb}
Let $\Xi$ be a a primitive proper $S$-adic subshift with unique decomposition property and endowed with the combinatorial metric. Then, it is $f$-embeddable.
\end{proposi}

\noindent It will be convenient to describe an $S$-adic subshift in terms of a Bratteli diagram (see for example \cite{DHS}). A weight on the Bratteli diagram permits to define an ultrametric on the set of its infinite paths. It will be necessary to prove that:
\begin{itemize}
  \item[--] the ultrametric Cantor set associated with this Bratteli diagram is bi-Lipschitz homeomorphic to the subshift with the combinatorial metric;
  \item[--] the weight on the Bratteli diagram satisfies the decay rate of Proposition~\ref{embed11.prop-criterion} so that the Cantor set associated with the diagram is embeddable.
\end{itemize}

\noindent The following result is needed.

\begin{lemma}[Durand~\cite{Dur1}, Lemma~8]
\label{embed11.lem-durand}
If the $S$-adic system generated by $(\sigma_n)_{n \in \NM}$ is
primitive with constant $s_0$, there exists a constant $K$ such that for all integers $r,s$, with $s-r \geq s_0$ and for all $b,c$ in $A_{s+1}$

$$\frac{|\sigma_{r,s+1}(b)|}{|\sigma_{r,s+1}(c)|} \leq K.
$$
\end{lemma}

\noindent Let $(\sigma_n: A_{n+1} \rightarrow A_n)$ be a primitive $S$-adic proper system with unique decomposition property and let $l$ and $r$ the letters associated with the properness. The Bratteli diagram is defined as follows

\begin{itemize}
  \item[--] For all $n \in \NM$, the set of vertices $\vs_n$ is in
bijection with the alphabet $A_n$: for each $a \in A_n$, there is a $v_a \in \vs_n$.

  \item[--] For all $n \in \NM$, an edge in $\es_n$ is a triple $e=(v_a,k,v_b)$ where $v_a \in \vs_n$, $v_b \in \vs_{n+1}$ and $k \in \NM$ such that the letter $a$ occurs in the word $\sigma_n(b)$ in position $k+1$. Hence, the number of edges from $v_a$ to $v_b$ is equal to the number of times the letter $a$ appears in $\sigma_n(b)$. 

  \item[--] If $e=(v_a,k,v_b)\in\es_n$ its {\em source} is $s(e)=v_a$, its {\em range} is $r(e)=v_b$ and its {\em label} is $\ell (e) = k$. 
\end{itemize}

\begin{exam}
\label{embed11.exam-edge}
{\em If $\sigma_n (b) = labcar$, then there are two edges $e_1,e_2$ from $v_a$ to $v_b$: one corresponds to the second letter and the other to the fifth letter of the word above. Their respective labels are $1$ and $4$.
}
\hfill $\Box$
\end{exam}
\begin{defini}
\label{embed11.def-path}
A path on the Bratteli diagram is a sequence of edges $(e_k)_{i \leq k < j}$ (with $1 \leq i < j  \leq +\infty$), such that $e_k \in \es_k$ for all $k$ and $r(e_k) = s(e_{k+1})$. Let $\Pi_n$ denote the set of paths with $i=1, \ j=n$ (simply called ``paths of length $n$''). Let $\Pi$ denote the union of all $\Pi_n$ and let $\Pi_\infty$ be the set of paths of infinite length ($i=1$, $j=+\infty$).
\end{defini}

\noindent Because the $\sigma_n$'s are taken from a finite set of substitutions $S$ the cardinality of the sets $\es_n$ and $\vs_n$ is uniformly bounded in $n$. Let

$$w_n = \big(\min\{|\sigma_{1,n}(a)| \ ; \ a \in A_n \}\big)^{-1}.
$$

\noindent It is a decreasing sequence which tends to $0$ as $n\to\infty$.
The weight of a finite path $\gamma$ will be defined by $w_{n+1}$ whenever $\gamma$ has length $n$. This leads to a metric on $\Pi_\infty$ defined by

$$d_w (x,y) = w_{n+1}, \quad \text{where } n 
 \text{ is the length of the longest common prefix of } x,y.
$$

\begin{proposi}
\label{embed11.prop-bilip}
The subshift $(\Xi,d)$, endowed with the combinatorial metric, and $(\Pi_\infty,d_w)$ are homeomorphic though a bi-Lipschitz homeomorphism.
\end{proposi}

\begin{proof} (i) {\bf Constructing a map $\Xi \rightarrow \Pi_\infty$.} Let $x = \ldots x_{-1} \cdot x_0 x_1 \ldots$ in~$\Xi$. Let $v_1 \in \vs_1$ be the vertex corresponding to the letter $x_0$. By the unique decomposition property, $x$ can be written in a unique way as the
concatenation

$$x =\ldots \sigma(x'_{-1})\underline{\sigma(x'_0)}\sigma(x'_1) \ldots\,,
$$

\noindent with $(x'_i)_{i \in \ZM} \in (A_2)^\ZM$ and such that the letter of $x$ of index $0$ is in the underlined word $\sigma(x'_0)$. Therefore, if $x'$ is the word $x' = (x'_i)_{i \in \ZM}$

\begin{equation}
\label{embed11.eq-psi}
 x = T^k \sigma (x'),
\end{equation}

\noindent for some $0 \leq k < | \sigma(x'_0) |$. This label $k$ corresponds to an occurrence of the letter $x_0$ in the word $\sigma(x'_0)$. It defines an edge $e=(v_1, k, v_2)$, where $v_2 \in \vs_2$ is the vertex corresponding to $x'_0$. Iterating this process of ``de-substitution'' leads to construct a sequence of words $x', x'', \ldots, x^{(n)} \ldots$ and a corresponding sequence
of edges $e_1, e_2, \ldots, e_n$. In particular $\psi(x) = (e_1, e_2, \ldots)$ defines a map $\Xi \rightarrow \Pi_\infty$. 

\vspace{.1cm}

\noindent (ii) {\bf $\psi$ is a bijection. } The map $\psi$ has an inverse $\phi: \Pi_\infty \rightarrow \Xi$ which will be built explicitly. Let $\gamma = (e_1, e_2, \ldots)$ be a path going through the vertices 
$v_1, v_2, \ldots$. Let $\epsilon$ be a symbol not in $A$ (the ``empty'' symbol). Now, given a finite word $w = (w_0, \ldots, w_{l-1}) \in A^\ast$, Define then $\bar{w} \in (A \cup \{\epsilon\})^\ZM$ as the sequence

$$\bar{w} = \ldots \epsilon \; \epsilon \; \epsilon \cdot w \; 
             \epsilon \; \epsilon \; \epsilon \ldots.
$$

\noindent This is a way of seeing a word in $A^\ast$ as a (partially defined) element of $A^\ZM$. For fixed $n$, define the sequence

\begin{equation}
\label{embed11.eq-phi}
\phi_n(\gamma) =
  T^{l(e_1)} \circ \sigma_1 \circ T^{\ell (e_2)} \circ 
   \ldots \circ \sigma_{n-2}\circ T^{\ell (e_{n-1})+1} 
    \left( \overline{r \; \sigma_{n-1} (a_n) \; l} \right)\,,
\end{equation}

\noindent where $\ell(e)$ denotes the label of the edge $e$. This sequence can be seen as an element of $(A \cup \{\epsilon\})^\ZM$, namely as:

$$\ldots \epsilon \epsilon \;
    \sigma_{1,n-1} (r) \; 
     \underline{\sigma_{1,n} (a_n)} \; 
      \sigma_{1,n-1} (l) \;
       \epsilon \epsilon \ldots,
$$

\noindent where the letter of index $0$ occurs in the underlined word. Its exact position is determined by the labels of the edges. This implies

\begin{equation}
\label{embed11.eq-phi2}
[\phi_n(\gamma)]_i \neq \epsilon 
  \quad \text{for} \quad 
  -|\sigma_{1,n-1}(r)| \leq i \leq |\sigma_{1,n-1}(l)|.
\end{equation}

\noindent By hypothesis, $\lim_{n \rightarrow +\infty} |\sigma_{1,n-1} (b)| = +\infty$ and this is true, in particular, whenever $b \in \{r,l\}$.
Therefore, the limit as $n$ tends to infinity of $\phi_n (\gamma)$ is an element of $A^\ZM$, which is noted $\phi(x)$. Using the definitions of $\phi$ and $\psi$ (see eq.~(\ref{embed11.eq-psi}) \&~(\ref{embed11.eq-phi})), it is straightforward that $\psi \circ \phi = \mathrm{id}_{\Pi_\infty}$. Therefore, $\psi$ is onto and $\phi$ is one-to-one. 

\noindent Conversely, let $x \in \Xi$, and $\gamma = \psi (x)$. It will be shown that $\phi(\gamma) = x$. For all $n$, $x$ has a unique decomposition

$$x = \ldots \sigma_{1,n} (b_{-1}) 
       \underline{\sigma_{1,n} (b_{0})} 
        \sigma_{1,n} (b_{1}) \ldots\,,
$$

\noindent with $b_i \in A_n$. So by definition of $\psi$ and $\phi$,

\begin{equation}
\label{embed11.eq-phin}
\phi_n (\gamma) = 
 \ldots \epsilon \; 
  \sigma_{1,n-1} (r) \; \underline{\sigma_{1,n} (b_0)} \; 
   \sigma_{1,n-1} (l) \; \epsilon \ldots\,,
\end{equation}

\noindent and the two words coincide on the word $\sigma_{1,n} (b_0)$ (it appears at the same position). Furthermore, by properness of the $S$-adic system, $\sigma_{n} (b_{-1})$ ends with the letter $r$ while $\sigma_{n}(b_1)$ begins with the letter $l$. Therefore, $\sigma_{1,n} (b_{-1})$ has $\sigma_{1,n-1} (r)$ as a suffix and $\sigma_{1,n} (b_{1})$ has $\sigma_{1,n-1} (l)$ as a prefix. So for all $i \in \ZM$ such that $[\phi_n(\gamma)]_i \neq \epsilon$, one has $[\phi_n(\gamma)]_i = x_i$. Taking a limit and using Equation~(\ref{embed11.eq-length}), $\phi(\gamma)$ and $x$ agree everywhere. Therefore, $\phi$ and $\psi$ are inverse bijections of each other.

\vspace{.1cm}

\noindent (iii) {\bf $\phi$ is bi-Lipschitz.} Let $\gamma, \gamma' \in \Pi_\infty$ be such that the $n$ first edges of $\gamma$ and $\gamma'$ coincide. Then, by definition, $\phi_n(\gamma) = \phi_n(\gamma')$. Thus $\phi(\gamma)$ and $\phi(\gamma')$ coincide for all indices $i$ satisfying $-|\sigma_{1,n-1}(r)| \leq i \leq |\sigma_{1,n-1}(l)|$. Using Lemma~\ref{embed11.lem-durand} and the definition of $w_n$ leads to

$$\forall b \in A_{n-1}\,,
\hspace{2cm}
   (w_{n-1})^{-1} \leq |\sigma_{1,n-1} (b)| \leq K (w_{n-1})^{-1}.
$$

\noindent Since $S$ is finite, $C = \max_{\sigma \in S, a \in A(\sigma)}|\sigma(a)|$ is well defined. Then, a word of the form $\sigma_{1,n} (b) = \sigma_{1,n-1} \circ \sigma_{n-1} (b)$ is at most $C$ times longer than the longest word of the form $\sigma_{1,n-1} (c)$. That is, using again Lemma~\ref{embed11.lem-durand}, $(w_n)^{-1} \leq CK(w_{n-1})^{-1}$.
If $[\phi(\gamma)]_i = [\phi(\gamma')]_i$ for $|i| \leq K (w_{n-1})^{-1}$, then in particular this inequality holds for $|i| \leq (w_n)^{-1}/C$. So

$$d(\gamma,\gamma') \leq w_n 
\hspace{1cm}\Rightarrow \hspace{1cm}
   d(\phi(\gamma), \phi(\gamma')) \leq Cw_n.
$$

\noindent Conversely, let $\gamma$ and $\gamma'$ coincide up to edge $n$, but differing on their $(n+1)$-th one. Then, by definition of $\phi$

\begin{eqnarray*}
\phi(\gamma) & =&
 \ldots \sigma_{1,n+1}(b_{-1}) 
  \underline{ \sigma_{1,n+1}(b_{0}) } \sigma_{1,n+1}(b_{1}) \ldots \\
\phi(\gamma') & = &
 \ldots \sigma_{1,n+1}(b'_{-1}) \underline{ \sigma_{1,n+1}(b'_{0}) }
  \sigma_{1,n+1}(b'_{1}) \ldots
\end{eqnarray*}

\noindent where $b_0 \neq b'_0$ and the letter of index $0$ belongs to the respective underlined words.
By the unique decomposition property, $\phi(\gamma)$ and $\phi(\gamma')$ have to differ at at an index $i$ satisfying

$$ -\max \{ |\sigma_{1,n+1}(b_{0})|, |\sigma_{1,n+1}(b'_{0})|\} 
   \leq i \leq  
    \max \{ |\sigma_{1,n+1}(b_{0})|, |\sigma_{1,n+1}(b'_{0})|\},
$$

\noindent that is for

$$ -K (w_n)^{-1} \leq i \leq K (w_n)^{-1}.
$$

\noindent This proves that $\phi^{-1}$ is Lipschitz.
\end{proof}

\vspace{.2cm}

\noindent {\bf Proof of Proposition~\ref{embed11.prop-bilip}. }In order to prove that $(\Xi,d)$ is $f$-embeddable, it is sufficient to
prove that $(\Pi_\infty,d_w)$ is embeddable. It suffices to remark that the corresponding Michon tree has the finite paths in $\Pi$ as vertices and the parent map is induced by deletion of the last edge. By finiteness of $S$, all paths in $\Pi_n$ have a bounded
number of extensions to paths of $\Pi_{n+1}$ and this bound is independent of $n$. Hence the Michon graph has a bounded number of children per vertex. The only thing left to show is that the sequence of weights $(w_n)_{n \in \NM}$ on the Bratteli diagram is bounded above by a decreasing geometric sequence. 

First, remark that for all $n$, the alphabet $A_n$ has at least two elements. Indeed, if one of the $A_n$ had only one element, the subshift $\Xi$ would be periodic (this case was ruled out). Using primitivity, for all letter $c \in A_{s_0+1}$, $|\sigma_{s_0}(c)| \geq 2$. By iteration, for all $k$ and all $c \in A_{ks_0 +1}$

$$|\sigma_{k s_0} (c)| \geq 2^k\,.
$$

\noindent So $w_{ks_0} \leq 2^{-k}$. Since $(w_n)_{n \in \NM}$ is decreasing,  $w_n \leq C \lambda^n$ with $\lambda = 2^{-1/s_0} < 1$, for some constant $C$. Thanks to Proposition~\ref{embed11.prop-criterion} it is proved that $(\Xi,d)$ has finite Assouad dimension, and is therefore embeddable.

 \subsection{An possible analogue result for tilings}
 \label{embed11.ssect-ana}

\noindent It is a natural question to ask whether linearly repetitive tilings of $\RM^n$ are $f$-embeddable. It seems reasonable to make the following conjecture (see for example~\cite{AC11} for the definitions of the objects). 
\begin{conj}
\label{embed11.conj-LRT}
Let $\Xi$ be the transversal of a tiling space of linearly repetitive tilings of $\RM^d$, together with the usual tiling metric. Then $\Xi$ is $f$-embeddable.
\end{conj}

\noindent The proof for subshifts is based on the fact that linearly repetitive subshifts have a good representation by Brattli diagrams.
While they are not self-similar, there are only a finite number of substitutions involved, which allows to generalize the methods of~\cite{JS11}. In a recent article, Aliste and Coronel~\cite{AC11} provide a good description of linearly repetitive tiling spaces of any dimension as the set of paths on a Bratteli diagram.
Given the transversal of a linearly repetitive tiling space, it is possible to represent it by the set of paths on a Bratteli diagram such that the number of vertices and edges in each $\vs_n$ or $\es_n$ is bounded uniformly in $n$. There are also have estimates for the size of the patches associated with paths of size $n$. These estimates are analogues of Durand's Lemma~\ref{embed11.lem-durand} and allow to define reasonable weights on the Bratteli diagram. In the end, the path space of this Bratteli diagram, endowed with the distance defined by the weight, is $f$-embeddable. However, there is no proof yet that the homeomorphism between this paths space and the tiling space is bi-Lipschitz.

\noindent In Proposition~\ref{embed11.prop-bilip}, the bi-Lipschitz character of $\phi$ is proved through using equation~\eqref{embed11.eq-phin}. This requires an estimate on the length of $\sigma_{1,n} (b_0)$, of $\sigma_{1,n-1} (r)$ and of $\sigma_{1,n-1} (l)$. This is a quantitative version of a property known as ``forcing the border''.
It seems that Aliste--Coronel's construction does satisfy a similar quantitative border forcing property, but this fact was not highlighted as such in their paper (as it seems to be a byproduct of their construction and not an essential feature of the result they prove). Establishing it would require a thorough reworking of their already sophisticated proof.
This justifies the fact that, even though this conjecture can be stated with some confidence, it is a conjecture and not as a result.

\section{Sturmian Sequences}
\label{embed11.sect-sturm}

 \subsection{Definitions and notations}
 \label{embed11.ssect-sturm-def}

\noindent First, a few facts need to be introduced about Sturmian sequences and their coding. This section follows in part the presentation in~\cite{Ar02}, and proof for some of the results cited below can be found there.

\begin{defini}
 Given a sequence $x=(x_n)_{n \in \ZM} \in A^{\ZM}$ on the finite alphabet $A$, define its language:
\[
 L_n(x) = \{ \text{finite words } a_1 \ldots a_n \text{ which appear in } x \}, \text{ and }
 L (x) = \bigcup_{n \in \NM} L_n (x).
\]
\end{defini}

\begin{defini}
 Given $x \in A^\ZM$, its complexity function $p_x$ is defined by
\[
 \forall n \in \NM, \quad p_x (n) = \mathrm{Card} (L_n).
\]
\end{defini}

\noindent When there is no risk of ambiguity, $p_x(n)$ is just noted $p(n)$. It is known that if $p_x(n) \leq n$ for some $n$, then the word $x$ is periodic.

\begin{defini}
A sequence is called \emph{Sturmian} if its complexity function satisfies $p(n) = n+1$ and if it is not eventually periodic, namely if it doesn't have a one-sided periodic infinite prefix or suffix).
\end{defini}

\noindent From the definition, a Sturmian sequence is a sequence on two letters (since $p(1) = 2$). These two letters are noted $0$ and $1$. Note that the ``not eventually periodic'' condition is here to rule out degenerate cases like $\ldots 111000 \ldots$ or $\ldots 0001000\ldots$.

\begin{proposi}
 The frequency of the letter $1$ in a Sturmian sequence $x$
\[
 \mathrm{freq}_x (1) := \lim_{n \rightarrow +\infty}\frac{\mathrm{Card}\{k \in [-n,n] \ ; \ x_k = 1\}}{2n+1}
\]
is well defined, and is an irrational number, noted $\alpha \in \left]0,1\right[$.
\end{proposi}

\begin{defini}
The Sturmian sequence $x$ is of \emph{type zero} if $\mathrm{freq}_x (1) < 1/2$, and of \emph{type one} if $\mathrm{freq}_x (1) > 1/2$.
\end{defini}

\noindent It is easy to see that in a Sturmian sequence of type $1$, the words $10$, $01$ and $11$ may appear, but not the word $00$. The same statement holds for sequences of type $0$, just exchanging the letters $0$ and $1$.

\vspace{.1cm}

\noindent Given a Sturmian sequence $x$, there is naturally a subshift of $\{0,1\}^{\ZM}$ associated with it.  It is by definition the set of all sequences $y$ such that $L (y) = L (x)$. It can also be defined as the closure of the orbit of $x$ in $A^\ZM$ (for the product topology). The two definitions are equivalent. It is of course shift-invariant (hence the name ``subshift''), and it is well known that it is minimal. If $x$ has frequency of ones equal to $\alpha$, so do all the elements of its subshift. In particular, it makes sense to write that a Sturmian subshift is of type $0$ or of type $1$. Conversely, the subshift associated with $x$ is exactly the set of all Sturmian sequences with the same frequency as $x$. Therefore, it makes sense to denote $\Xi(\alpha)$ the subshift associated with $\alpha \in (0,1)$. Sturmian sequences can be \emph{recoded} using the following substitutions. Define:
\[
 \sigma_0: \left\{ \begin{array}{ll} 0 & \mapsto 0 \\ 1 & \mapsto 10 \end{array} \right. \quad
\text{and} \quad 
 \sigma_1:
\left\{ \begin{array}{ll} 0 & \mapsto 01 \\ 1 & \mapsto 1. \end{array} \right.
\]

\begin{proposi}
\label{embed11.prop-recoding}
 For any Sturmian sequence $x$ of type $0$, there is a Sturmian sequence $x'$ such that either $x = \sigma_0(x')$ or $x = T \sigma_0 (x')$, where $T$ is the shift operator.

\noindent For any Sturmian sequence $x$ of type $1$, there is a Sturmian sequence $x'$
such that either $x = \sigma_1(x')$ or $x = T \sigma_1 (x')$.
\end{proposi}

\noindent Let $\Phi$ denote the recoding map $x \mapsto x'$.  It is worth noting that if $x,y$ are two elements of the same Sturmian subshift, then $\Phi(x)$ and $\Phi(y)$ belong to the same subshift. 
Non-periodicity of Sturmian sequences implies that for any Sturmian sequence $y$ of type $0$ (resp.\ $1$), there is a $k$ such that $\Phi^k(y)$ is of type $1$ (resp.\ $0$).

\begin{defini}
 Let $x$ be a Sturmian sequence, and $(\Phi^n(x))_{n \in \NM}$ be the sequence of recoded Sturmian sequences. By definition, for all $n$,
\[
 \sigma_0^{b_1} \circ \sigma_1^{b_2} \circ \ldots \circ \sigma_0^{b_{2n+1}} (\Phi^{2n+1} (x))
\text{ and }
\sigma_0^{b_1} \circ \sigma_1^{b_2} \circ \ldots \circ \sigma_1^{b_{2n}} (\Phi^{2n} (x))
\]
are in the same orbit as $x$. All the $b_n$ are positive, except maybe $b_0 = 0$. The sequence $(b_n)_{n \in \NM}$ is called the \emph{multiplicative coding} of the Sturmian sequence $x$.
\end{defini}

\noindent All Sturmian sequences in a same subshift have the same multiplicative coding, and an acceptable multiplicative coding determines uniquely a Sturmian subshift.

 \subsection{Partial fraction decomposition and multiplicative coding}
 \label{embed11.ssect-partialfrac}

\noindent The properties of a Sturmian subshift are closely related to the partial fraction decomposition of the number $\alpha$ (which is the frequency of ones in the subshift).

\noindent Let $\alpha \in \RM$ be irrational. Then $\alpha$ can be written uniquely $\alpha = a_0 + \alpha_0$, with $a_0 \in \ZM$ and $\alpha_0 \in (0,1)$. The Gauss map $G:[0,1]\to[0,1]$, applied to $\alpha_0$, generates the continuous fraction expansion

$$G(\alpha) = \frac{1}{\alpha} -a(\alpha)\,,
\hspace{2cm}
   a(\alpha) = \left[\frac{1}{\alpha}\right]\,,
$$

\noindent where $[x]$ denotes the {\em integer part} of $x$, namely the largest integer smaller than or equal to $x$. Hence

$$\alpha= a_0 +
   \frac{1}{a_1+\alpha_1}\,,
\hspace{2cm}
   a_1=a(\alpha)\,,\;\; \alpha_1=G(\alpha)\,.
$$

\noindent Iterating this formula gives rise to the continuous fraction expansion

$$\alpha = a_0 +
   \cfrac{1}{a_1+
    \cfrac{1}{a_2+
     \cfrac{1}{a_3+\dotsb+
      \cfrac{1}{a_n+\alpha_n
}}}}\,,
\hspace{2cm}
  \alpha_{n+1}=G(\alpha_n)\,,\; a_n=a(\alpha_{n-1})\,\;\;n\geq 1\,.
$$

\noindent The standard notation is 

$$\alpha = [a_0;a_1,a_2,\cdots, a_n,\cdots]\,,
$$

\noindent and the $a_n$'s are called the {\em partial quotients} of $\alpha$.

\begin{defini}
\label{embed11.def-bndtype}
A number $\alpha$ has {\em bounded type} whenever, the sequence of its partial quotients is bounded.
\end{defini}

\noindent Having bounded type is an exceptional property. This is a theorem of Khintchine~\cite{Kh97}. See also Levy~\cite{Le36}.

\begin{theo}
\label{embed11.th-contfrac}
For almost every $\alpha\in[0,1]$ the sequence $(a_n)_{n\in\NM_\ast}$ of partial quotients of $\alpha$ is unbounded.
\end{theo}

\noindent One famous example of such a typical number is $e=2,71828\cdots$ the continued fraction of which was computed by Leonhard Euler in 1737 namely

$$e=[2,1,2,1,1,4,1,1,6,\cdots, 1,1,2l,\cdots]
\hspace{2cm}
  a_{3l-1}=2l\,,\;\; a_{3l-2}=a_{3l}=1\,,\;\;l\geq 1\,.
$$

\noindent Many properties of a Sturmian subshift are determined by the arithmetic properties of the number $\alpha$ associated with it. The following theorem is proved in Hedlund and Morse's seminal paper.

\begin{theo}[see \cite{MH40}]
\label{embed11.th-sturm-LR}
Consider a Sturmian subshift $\Xi \subset \{0,1\}^\ZM$, of parameter $\alpha$.
This subshift is \emph{linearly repetitive} (see Definition~\ref{embed11.def-LR}) if and only if $\alpha$ has bounded type.
\end{theo}

\noindent An immediate consequence of this theorem and of Theorem~\ref{embed11.th-LR} is the following result.

\begin{proposi}
\label{embed11.prop-sturm-embed}
Let $\Xi$ be a Sturmian subshift with parameter $\alpha$, endowed with the combinatorial metric. If $\alpha$ has bounded type, then $\Xi$ is $f$-embeddable.
\end{proposi}

\noindent This result needs a converse statement: if $\alpha$ has not bounded type, then the associated Sturmian subshift is not $f$-embeddable. This will be proved in Section~\ref{embed11.ssect-sturm-nonembed}.

\noindent One way to understand the deep links between the combinatoric properties of a Sturmian subshift and the arithmetic properties of its parameter $\alpha$ is the following. A Sturmian sequence $x$ can be seen as the coding of an orbit of the rotation of angle $\alpha$ on the circle: there is $\{I_0,I_1\}$ a partition of the circle, and $s$ a point in the circle such that $x_n = 0$ if and only if $(s+\alpha \mod 1) \in I_0$.
The rotation on the circle can be related to the continued fraction decomposition of $\alpha$ on the one hand, and on the multiplicative coding of $x$ on the other hand, to get to following result.

\begin{theo}[(see \cite{Ar02}]
Let $x$ be a Sturmian sequence, and $\alpha = [0;a_1,a_2, \ldots]$ be the frequency of $1$ in $x$.
Then the multiplicative coding of $x$ is given by the partial quotients of $(1-\alpha) / \alpha = \alpha^{-1} - 1$.
\end{theo}

\noindent It is straightforward that $(1-\alpha)/\alpha = [a_1-1; a_2, a_3, \ldots]$. In particular, $\alpha$ has bounded type if and only if the sequence of coefficients of the multiplicative coding is bounded.

 \subsection{Non embeddability of certain Sturmian subshifts}
 \label{embed11.ssect-sturm-nonembed}

\noindent This section is devoted to the proof of the following result.

\begin{proposi}
\label{embed11.prop-sturm-nonembed}
 Consider a Sturmian sequence $x$, with associated subshift $\Xi$ and associated parameter $\alpha$. If $\alpha$ has unbounded partial quotients, then $\Xi$ (with the combinatorial metric) is not embeddable in a finite dimensional space.
\end{proposi}

\noindent A Sturmian subshift $\Xi$ is well described combinatorially by the (bilateral) tree of words of its elements.

\begin{defini}
 Given a Sturmian subshift $\Xi$, its un-reduced tree of words is defined as follows.
\begin{itemize}
 \item[--] For all $n \geq 0$, define a sequence of refining partitions of $\Xi$ by:
\[
 \Pp_n = \{ [y_{-n} \ldots y_n ] \ ; \ y \in \Xi \},
\]
where $[y_{-n} \ldots y_n ]$ is the cylinder set of all words in $\Xi$ which coincide with $y$ on indices $-n \leq i \leq n$. By convention, $\Pp_{-1}:= \{\Xi \}$.
 \item[--] The set of vertices of the tree is in bijection with the disjoint union of all the $\Pp_n$. If $U \in \Pp_n$, the associated vertex is noted $v_U$.
 \item[--] The ancestor relation is induced by inclusion: $v_U \preceq v_V$ if and only iv $U \subseteq V$.
 \item[--] The weight of the vertex $v_U$ is $1/(n+2)$ if $X \in \Pp_n$ ($n \geq -1$).
\end{itemize}
The (reduced) tree of words is obtained from this tree by the reduction process defined in section~\ref{embed11-ssect-trees}.
\end{defini}

\noindent One remark about this tree: if $X \in \Pp_n$, then $\diam(X) \leq n^{-1}$, for the combinatorial distance, with equality if and only if $X$ is the non-trivial union of two distinct elements of $\Pp_{n+1}$. Note that the vertices $v_X$ of such clopen sets $X$ are exactly the vertices with two children: these are precisely the ones which are not dropped by the reduction process.

\noindent This leads to the following proposition.

\begin{proposi}
 The boundary of the reduced tree of words associated with a Sturmian subshift $\Xi$ is bi-Lipschitz homeomorphic to the subshift (endowed with the combinatorial metric).
\end{proposi}

\begin{proof} The proof is almost tautological: given an infinite path in the tree, say $\gamma$, the sequence of vertices $(v_{X_{n}})_{n \in \NM}$ defines a decreasing sequence of compact sets $X_n$, the diameter of which tends to zero. Therefore, its intersection is not empty and consists of a single element $\{x\}$. Define $\phi(\gamma) = x$.

\noindent Conversely, given $x \in \Xi$, there is a unique decreasing sequence of sets $X_n \in \Pp_n$, such that for all $n$, $x \in X_n$ (explicitly: $X_n = [x_{-n} \ldots x_n]$). Then it defines an infinite path $\gamma \in \partial \ts$. Clearly, $\gamma$ is the unique pre-image of $x$ by $\phi$.

\noindent The fact that $\phi$ is bi-Lipschitz (and in particular, it is a homeomorphism) results from the remark above, on the diameters of the elements of $\Pp_n$.
\end{proof}

\vspace{.2cm}

\begin{lemma}
 Let $y$ be a Sturmian sequence, assume that $y = \sigma_0^{b_n}(z)$, with $z$ a Sturmian sequence. Then $y$ contains the words $10^{b_n}1$ and $10^{b_n}0$.
\end{lemma}

\begin{proof} If $y = \sigma_0^{b_n} (z)$, then $y$ is of type $0$ and $z$ is of type $1$. By iteration, it is straightforward that $\sigma_0^{b_n}(0) = 0$ and $\sigma_0^{b_n} (1) = 1 0^{b_n}$. Since $z$ is of type $1$, it contains the words $10$ and $11$. Therefore, $y$ contains the words:
\[
 10^{b_n}10^{b_n} \quad \text{and} \quad 10^{b_n}0.
\]

\end{proof}

\vspace{.2cm}

\noindent {\bf Proof of Proposition~\ref{embed11.prop-sturm-nonembed}: } Let $x$ be a Sturmian sequence, with $\alpha$ the frequency of $1$, and assume that $\alpha$ does not have bounded type. Then its multiplicative coding $(b_n)_{n \in \NM}$ is an unbounded sequence.

\noindent Let $y$ = $\Phi^{b_0+b_1+\ldots+b_{n-1}} (x)$, and $z = \Phi^{b_n}(y)$. Without loss of generality, we assume that $z$ is of type $1$, so that
\[
 y = \sigma_0^{b_n}(z).
\]
\noindent Then, using previous lemma, $y$ contains the words $10^{b_n}1$ and $10^{b_n}0$. Therefore, for all $1 \leq k \leq [b_n]/2 -1$, $y$ contains the words $0 0^{2k} 0$ and $0 0^{2k} 1$.

\noindent Applying the substitution $\sigma := \sigma_0^{b_1} \circ \ldots \circ \sigma_0^{b_{n-1}}$, the sequence $x$ contains the words
\[
 \sigma(0)\sigma(0)^{2k} \sigma(0) \quad \text{and} \quad \sigma(0)\sigma(0)^{2k} \sigma(1).
\]
\noindent Let $a$ be the last letter of $\sigma(0)$ and $\sigma(1)$ (it is the same: $0$ if $b_0 \neq 0$, $1$ otherwise). Since $\sigma(1)$ starts by $1$ and $\sigma(0)$ starts by $0$, $x$ contains the words
\begin{equation}\label{eq-arbre-sturmien}
 a\sigma(0)^k\sigma(0)^k0 \quad \text{and} \quad a\sigma(0)^k\sigma(0)^k1.
\end{equation}
\noindent Then, let $X_k = [a \sigma(0)^k \cdot \sigma(0)^k]$, where the dot separates the indices $i \leq 0$ and $i > 0$. It is an element of $\Pp_{k |\sigma(0)|}$, where $|\sigma(0)|$ is the length of $\sigma(0)$. Let $v_k$ be the associated vertex. Note that $v_{k+1}$ is a child of $v_k$ for all $k$.
From Equation~\eqref{eq-arbre-sturmien}, the vertices $v_k$ have two distinct children (in the non-reduced tree of words), therefore, they are elements of the reduced tree, and their weight is $(k |\sigma(0)|)^{-1}$. In particular, this construction shows that there are two vertices $u,v$ (namely $v_1$ and $v_{[b_n]/2-1}$), such that the quotient of their weights is $[b_n/2]-1$, and their distance in the tree is $[b_n/2]-1$.

\noindent This construction can be done for all $n$. If $(b_n)_{n \in \NM}$ is unbounded, this shows that the weights cannot satisfy the geometric decay condition of Proposition~\ref{embed11.prop-criterion}, and $\Xi$ is not embeddable.
\hfill $\Box$

\vspace{.2cm}

\noindent The results presented here (namely proposition~\ref{embed11.prop-sturm-embed} and \ref{embed11.prop-sturm-nonembed}) provide a proof of Theorem~\ref{embed11.th-sturmHull}. It is actually possible to give a more complete version of it.

\begin{theo}
\label{embed11.th-}
A Sturmian subshift $\Xi(\alpha)$ is $f$-embeddable if and only if the irrational number $\alpha$ associated with it has bounded type. In particular
\begin{enumerate}
 \item[(i)] if $\alpha$ is a quadratic irrational, then $\Xi(\alpha)$ is $f$-embeddable;
 \item[(ii)] for almost every $\alpha \in (0,1)$, the subshift $\Xi(\alpha)$ is not $f$-embeddable;
 \item[(iii)] the boundary of $\Xi(e)$ is not embeddable for $e=2.71828\ldots$.
\end{enumerate}
\end{theo}

\begin{proof}
The first part of this result is Theorem~\ref{embed11.th-sturmHull} which was proved above. Point (i) is a consequence of the fact that quadratic irrational have an eventually periodic (hence bounded) continued fraction expansion.
Point (ii) is a consequence of the theorem of Khintchine on continued fractions, and point (iii) is a consequence of the explicit formula for the continued fraction decomposition of $e$.
\end{proof}

\section{Hausdorff Dimension}
\label{embed11.sect-haus}

\noindent The Assouad dimension provides a good measure of embeddability of ultrametric Cantor spaces in $\RM^n$. In~\cite{JS11}, it was proved that for \emph{self-similar} Cantor sets, the smaller dimension of a space $\RM^n$ in which an ultrametric Cantor $(C,d)$ set can be embedded is equal to $\lfloor \dim_\mathrm{H} (C,d) \rfloor +1$, where $\dim_\mathrm{H} $ is the Hausdorff dimension.
When the space is not self-similar however, there is no reason that the Hausdorff dimension has anything to do with the possibility to embedd the space. For example~\cite{LM94} mentions an example of a countable metric space (hence of Hausdorff dimension $0$) which cannot be bi-Lipschitz embedded in any finite-dimensional Euclidean space. The goal of this section is to revisit the definition and properties of Hausdorff dimension in terms of the Michon representation for ultrametric Cantor sets. Proposition~\ref{embed11.th-HausCpart} gives an adapted formula for the computation of the Hausdorff dimension, which will be used in next section.

\subsection{Hausdorff Dimension of the boundary of a tree}
\label{embed11.ssect-comp}

\noindent The definition of the Hausdorff dimension~\cite{Fa90} starts with the following construction: given an open cover $\us$ of $C=\partial\ts$ and given $s\in [0,\infty)$, let $\Hh^s(\us)$ be defined by

\begin{equation}\label{embed11.eq-hausdorff-cover}
\Hh^s(\us) =
   \sum_{U\in\us} \diam(U)^s\,.
\end{equation}

\noindent In addition, let $\diam(\us)=\sup_{U\in\us} \diam(U)$. Then for $0<\delta<1$ let $\Hh_\delta^s(C)$ be defined by

\begin{equation}
\label{embed11.eq-hdelta}
\Hh_\delta^s(C) =
   \inf_{\idiam(\us)<\delta} \Hh^s(\us)
\end{equation}

\noindent It follows that $\delta'\leq\delta\Rightarrow \Hh_\delta^s(C)\leq \Hh_{\delta'}^s(C)$. In addition, if $\sigma >0$ then $\Hh_\delta^{s+\sigma}(C)\leq \delta^\sigma\Hh_\delta^s(C)\leq \Hh_\delta^s(C)$. Consequently, 

(i) the limit $\lim_{\delta\to 0} \Hh_\delta^s(C)=\Hh^s(C)$ exists in $[0,+\infty)\cup\{+\infty\}$,

(ii) there is a unique $s_0$ such that if $s>s_0$ then $\Hh^s(C)=0$ whereas for $s<s_0$, $\Hh^s(C)=+\infty$.

\noindent This unique value $s_0$ is precisely the Hausdorff dimension $s_0=\dim_\mathrm{H} (C,d)$. In order to give a more tractable formula, the following definition will be nedded

\begin{defini}
\label{embed11.def-fintree}
A finite subtree $\Gamma$ of $\ts$ is a tree graph $\Gamma=(\vs_\Gamma, \es_\Gamma, \bullet)$ where
\begin{enumerate}
 \item[(i)] $\vs_\Gamma\subset \vs$ is finite and contains the root $\bullet$;
 \item[(ii)] for every $v\in\vs_\Gamma$ every ancestor of $v$ belong to $\vs_\Gamma$;
 \item[(iii)] the ``parent'' relation is induced by the one on $\ts$ (equivalently, the edges of $\Gamma$ are the edges of $\ts$ between vertices of $\Gamma$).
\end{enumerate}
$\Ll \Gamma$ will denote the set of leaves (vertices of maximal depth).
The finite subtree $\Gamma$ is called full if for any vertex $v\in\Gamma \setminus \Ll \Gamma$, each child of $v$ is a vertex of $\Gamma$. It is equivalent to the statement $\partial \ts = \bigcup_{v \in \Ll \Gamma}[v]$.
\end{defini}

\begin{figure}[ht]
   \centering
\includegraphics[width=9cm]{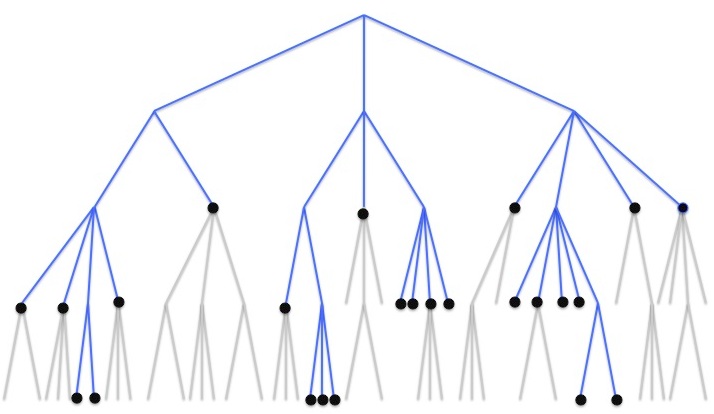}
\caption{A full finite subtree and its extremal vertices}
\label{comp11.fig-fulltree}
\end{figure}

\noindent The following theorem allows for simpler computation of the Hausdorff dimension.

\begin{proposi}
\label{embed11.th-HausCpart}
Let $\ts$ be a reduced Michon tree with weight $\kappa$. For any $\delta >0$ let $\GG_\delta(\ts)$ be the set of all full finite subtrees of $\ts$ such that $\max_{v\in\Ll \Gamma}\kappa(v) <\delta$. Then
\begin{equation}\label{embed11.eq-hausdorff-subtree}
\Hh^s_\delta(\partial\ts) =
   \inf_{\Gamma\in\GG_\delta(\ts)} 
    \sum_{v\in\Ll \Gamma} 
     \kappa(v)^s\,.
\end{equation}
\end{proposi}

\noindent The following theorem was known, even for more general metric spaces than ultrametric ones. With our notations, the proof becomes very simple and is included.
\begin{coro}[Assouad~\cite{As79}, Prop.\ 2-(j)]
\label{embed11.cor-upbnd}
Let $(X,d)$ be an ultrametric Cantor set.
Then
\[
 \dim_\mathrm{H} (X,d) \leq \dim_\mathrm{A} (X,d).
\]
\end{coro}

\noindent Conversely, for any reduced tree $\ts$, there exists a weight which turns $\partial \ts$ into an ultrametric Cantor set of Hausdorff dimension $1$. It is the Kraft weight defined inductively by $\kappa (\bullet)=1$ and for all $w$, if one notes $v$ the parent of $w$ and assume that $v$ has $k$ children, then $\kappa (w) = \kappa (v) /k$. It leads to the following corollary, the proof of which is left to the reader,

\begin{coro}
\label{embed11.cor-kraftw}
Let $\ts$ be a reduced Michon tree with its Kraft weight $K$.
The Hausdorff dimension of $(\partial \ts, d_K)$ is exactly one.
In particular $(\partial \ts, d_K)$ is $f$-embeddable if and only if $\ts$ has a bounded number of children per vertex.
\end{coro}

\noindent The proofs of the first two results will be the content of the following section.

\subsection{Proof of Proposition~\ref{embed11.th-HausCpart}}
\label{embed11.ssect-proof}

\noindent Remark that the right-hand side of formula~\eqref{embed11.eq-hausdorff-subtree} is an infimum of the quantity $\sum_{u \in \Uu}\diam(U)^s$ taken over the special partitions of the form $\Uu=\{[v] \ ; \ v \in \Ll \Gamma\}$. Therefore,
\begin{equation}\label{embed11.eq-hausdorff-ineq}
 \Hh_\delta^s(\partial \ts) \leq \inf_{\Gamma\in\GG_\delta(\ts)} \sum_{v\in\Ll \Gamma} \kappa(v)^s\,.
\end{equation}

\noindent Conversely, let $\Uu$ be a cover of $\partial \ts$ by open sets of diameter less than $\delta$. By definition of the distance on $\partial \ts$, for each $U \in \Uu$, there is a vertex $v_U$ such that $U \subset [v_U]$ and $\kappa(v_U) = \diam (U)$.
The family $(v_U)_{U \in \Uu}$ form a covering of $\partial \ts$ by clopen sets of diameter less than $\delta$. By compactness, there is a finite sub-cover $[v_1], \ldots, [v_m]$. Then, for each pair $v_i \neq v_j$, exactly one of the following holds: (a) $v_i$ is an ancestor of $v_j$ (in which case $[v_j] \subset [v_i]$); (b) or the converse holds; (c) or neither is an ancestor of the other (in which case $[v_i] \cap [v_j] = \emptyset$). By removing all clopen sets in this covering which are included in others, it is possible to get a finite \emph{partition} by clopen sets $([v'_1], \ldots, [v'_n])$. In the process, terms have been removed. Therefore
\[
 \sum_{U\in\Uu} \diam(U)^s = \sum_{U\in\Uu} \kappa(v_U)^s \geq \sum_{i=1}^n \kappa(v'_i)^s.
\]
It is straightforward to see that the subtree consisting of vertices $v'_1, \ldots, v'_n$ and all their ancestors is a finite full subtree of $\ts$.
Therefore, for each covering of $\partial \ts$ by open sets, there is a full finite subtree $\Gamma$ such that
\[
\sum_{U\in\Uu} \diam(U)^s \geq \sum_{v\in\Ll \Gamma} \kappa(v)^s.
\]
It proves that inequality~\eqref{embed11.eq-hausdorff-ineq} is in fact an equality, which achieves the proof.

\subsection{Proof of Corollary~\ref{embed11.cor-upbnd}}

\noindent Let $(X,d)$ be an ultrametric Cantor set represented by its reduced Michon tree $\ts$ with weight $\kappa$. Without loss of generality, the weight of the root $\kappa(\delta)$ can be taken to be $1$. If $(X,d)$ is $s$-subhomogeneous, then, for any $\varepsilon > 0$, $\Hh^{s+\varepsilon} (X) = 0$. For indeed, using Proposition~\ref{embed11.prop-shomMichon}, there is $C$ such that for all $\delta$, the following holds (where $\bullet$ is the root of $\ts$):
\[
  \# \Big( \Ll (\ts,\bullet,\delta) \Big) \leq C \delta^{-s}.
\]

\noindent It ought to be remarked that $\Ll (\ts,\bullet,\delta)$ is nothing but the leaves of the full finite subtree $\Gamma$ consisting precisely of these leaves and all of their ancestors up to the root. By definition, these leaves define clopen sets $[v]$ of diameter less than $\delta$. Therefore by Proposition~\ref{embed11.th-HausCpart},
\[
 \Hh_\delta^{s+\varepsilon} (X) \leq \sum_{v \in \Ll\Gamma} \kappa(v)^{s+\varepsilon}.
\]
Now, the number of terms in the sum is bounded above by $C \delta^{-s}$. Therefore $\Hh_\delta^{s+\varepsilon}$ is bounded above by $C\delta^{\varepsilon}$, which tends to $0$ as $\delta$ tends to~$0$.

\section{Random Trees}
\label{embed11.sect-GW}

\noindent Random trees are an ubiquitous objects in modern mathematics. It is therefore natural to investigate the boundary of such trees, at they are good candidates to provide examples of Cantor sets with typical properties. Originally, Sir Francis Galton and the reverend Watson introduced these trees to investigate the probability of disappearance of family names among the British aristocrats in the nineteen century.
Nowadays, Galton-Watson processes appear in situations as varied as, for example, the description of the nuclear chain reactions or of electron emission in a photomultiplier tube (see~\cite{Ha63}). 

\vspace{.1cm}

\noindent A Galton-Watson branching process, as it is called today~\cite{Ha63,AN72} defines naturally a tree. Starting with a root $\bullet$, its number of children $\xi_\bullet$ is a random variable with integer value and probability distribution $p$. Each offspring is represented by a vertex $v$ and each vertex has a random variable $\xi_v$ attached to it, describing the number of its children. All $\xi_v$ are independent and identically distributed, with distribution $p$. Going inductively in the same way, leads to a tree. It is finite with positive probability if $p_0:=p(\{0\}) > 0$ and is almost surely infinite otherwise. In order to fit in the setting of this paper, the resulting tree ought to be infinite and reduced namely each vertex should have at least two children. In order to ensure it, it will be assumed that $p_0 = p_1 = 0$. The resulting branching process is called a \emph{reduced random tree}.

\vspace{.1cm}

\noindent It is worth noting that, using the deterministic weight $\kappa(v) = 2^{- |v|}$, the dimension of the boundary of trees obtained by a Galton-Watson process was computed by Hawkes~\cite{Haw81} to be given by $\log(m)/ \log(2)$, where $m > 1$ is the average number of children per vertex. The Hausdorff dimension associated with this metric can be used to define an ``average branching number''-even for more general trees- which turns out to be the crucial parameter when studying the behavior of random walks on trees~\cite{Lyo92}. The approach taken in the present paper is however different, since the trees are endowed with random weights rather than a deterministic one. To the authors best knowledge, this approach is new.

 \subsection{Reduced Random Trees: a review}
 \label{embed11.ssect-GWrev}

\noindent Let $p=(p_n)_{n\geq 2}$ be the probability supported by $[2,\infty)\subset \NM$ describing the number of children of each vertex. Let $Z_n$ be the number of descendant at the generation $n$. Following Watson's idea, it is convenient to introduce the generating function 

$$P_n(x) =
  \sum_{l=2}^\infty \pr\{Z_n=l\}\;x^l
$$

\noindent The construction of the random tree implies the following formula, where $\vs_n$ denotes the set of vertices at generation $n$,

\begin{equation}
\label{embed11.eq-GWind}
Z_{n+1}=\sum_{v\in\vs_n} \xi_v\,.
\end{equation}

\noindent Consequently the conditional probabilities are given by

$$\pr\{Z_{n+1}=l\,|\,Z_n=k\}=
   \sum_{j_1+\cdots+j_k=l} p_{j_1}\cdots p_{j_k}\,.
$$

\noindent In particular, it shows that $(Z_n)_{n\in\NM_\ast}$ defines a Markov chain. Moreover

\begin{equation}
\label{embed11.eq-GWrec}
P_{n+1}(x) = P_n\big(P(x)\big)\,,
\hspace{2cm}
  P(x) = \sum_{n=2}^\infty p_n\,x^n\,.
\end{equation}

\noindent Since $p$ is a probability, the series defining $P$ converges for $0\leq x\leq 1$. In addition, $P(1)=1$, $m=P'(1)=\EM(\xi)$ represents the average number of offsprings, if it exists. With the present restrictions, it follows that $m\geq 2$. The function $x\in[0,1]\to  P(x)\in[0,1]$ is positive, increasing and convex, more generally all its derivatives, when they exist, are positive, meaning that $P$ is {\em completely monotone}. From the recursion relation~(\ref{embed11.eq-GWrec}) it follows immediately that

$$P_n(x) = \underbrace{P\circ P\circ\cdots\circ P}_{n}(x)
$$

\noindent In particular, since $P(1)=1$, it follows that $\EM(Z_n)= P_n'(1)=m^n$ for all $n$. Hence the number of descendants at generation $n$ grows exponentially fast with the generation in the average.

\vspace{.1cm}

\noindent Let now $\Ff_n$ be the sigma algebra generated by the variables $\xi_v$ for $v$ vertices of the generations $ k\leq n$. A classical remark made by Doob \cite{Do53}, is that, thanks to the equation~(\ref{embed11.eq-GWind}) defining the process, the family $(Z_n)_{n\in\NM_\ast}$ satisfies

$$\EM\left(Z_{n+1}|\Ff_n\right)= mZ_n\,.
$$

\noindent Namely $W_n=Z_n/m^n$ is a {\em martingale}. As shown earlier by Doob \cite{Do40} this implies

\begin{theo}[See \cite{Ha63,AN72}]
\label{embed11.th-martingale}
Let $p=(p_n)_{n\geq 2}$ be the probability distribution for the number $\xi$ of offsprings such that the average $\EM(\xi)=m$ and the variance $\var(\xi)=\sigma^2$ are finite. Then the sequence $W_n=Z_n/m^n$ of random variables is a martingale with respect to the increasing sequence $\Ff_n$ of $\sigma$-algebras. In particular it converges almost surely to a random variable $W$ such that

$$\EM(W)=1\,,
\hspace{2cm}
  \var(W) = \frac{\sigma^2}{m^2-m}\,.
$$
\end{theo}

 \subsection{Proof of the Proposition~\ref{embed11.prop-GaltWatson}}
 \label{embed11.ssect-propGW}

\noindent If the probability $p$ has an infinite support, given any integer $M\geq 2$, the probability $P_M$ that a given vertex has more than $M$ children is non zero. The construction of the Random Reduced Tree, can be seen by associating inductively with any vertex of generation $n$ a string $(b_1,b_2,\cdots, b_n)$ of integers so that $1\leq b_n\leq \xi_{b_1,\cdots, b_{n-1}}$. In particular, since $\xi_v\geq 2$ for all $v$'s almost surely the subset $\ws\subset\vs$ made of vertices for which all $b_j$'s belong to $\{1,2\}$ is non empty and gives an infinite binary subtree. Since the random variables $\{\xi_v\,;\, v\in\ws\}$ are {\em i.i.d.}, it follows that , given $M\geq 2$ 

$$\pr\{\xi_v\leq M\,;\, \forall v\in A\}=
   (1-P_M)^{\# A}\,.
$$

\noindent In particular the probability that all vertices in $\ws$ have $\xi_v\leq M$ vanishes.
\hfill $\Box$

 \subsection{Random Weight}
 \label{embed11.ssect-rweight}

\noindent In order that the boundary of the rooted random tree $\ts=(\vs,\es,\bullet)$ previously built becomes an ultrametric Cantor set, it is necessary to put a weight on each vertex. The previous construction suggests that the weight itself be random and Markovian as well. In order to do so, the following model of random weight is proposed: let $(\lambda_v)_{v\in\vs}$ be a family of {\em i.i.d} in $[0,1]$ with common distribution $\rho$; then the weight $\kappa(v)$ will be given by

$$\kappa(v)=\lambda_v\,\kappa(u)\,,
\hspace{1cm}\mbox{\rm if}\;
   u=\mbox{\rm parent of }\; v\,,
\hspace{2cm}
    \kappa(\bullet)=1\,.
$$

\noindent In order that this defines a good weight, it is required that $\lambda_v\neq 0$ with probability one, namely $\rho\{0\}=0$. Moreover, in order that the weight converges to zero along any infinite paths with probability one, it will be required that $\rho\{1\}<1$. It will be convenient to use the following generating function (Mellin transform)

$$h(s)=
  \int_0^1\lambda^s \rho(d\lambda)\,.
$$

\noindent It is easy to see that $h$ is decreasing, logarithmically convex, namely $h(\eta s_1+(1-\eta)s_0)\leq h(s_1)^\eta\,h(s_0)^{1-\eta}$ for $0<\eta <1$, and $\lim_{s\to\infty}h(s)=\rho\{1\}$.

 \subsection{Proof of Theorem~\ref{embed11.th-GWHaus}}
 \label{embed11.ssect-theoGWHaus}

\noindent Thanks to Proposition~\ref{embed11.th-HausCpart}, the Hausdorff dimension of the tree can be computed from the following random variables

\begin{equation}
\label{embed11.eq-hausmeas}
\Hh^s(\Gamma) =
   \sum_{v\in\Ll\Gamma} \kappa(v)^s\,,
\hspace{2cm}
   \Hh_n^s= \sum_{v\in\vs_n} \kappa(v)^s\,.
\end{equation}

\noindent where $\Gamma$ is a full finite subtree of the random $\ts$. In the following $\GG$ will denote the set of finite full subtrees of $\ts$. This set is ordered by the inclusion of the vertex sets. In addition, giving $\Gamma,\Gamma'\in\GG$, then they admit a least upper bound $\Gamma\vee\Gamma'$ and a greatest lower bound $\Gamma\wedge \Gamma'$, as can be checked immediately. This suggests to define a new family of $\sigma$-algebras: let $\hFf_\Gamma$ be the $\sigma$-algebra generated by the $\xi_u$'s and the $\lambda_v$'s where $v\in\vs_\Gamma$ and $u\in\vs_\Gamma\setminus\Ll\Gamma$. In particular, if $\Gamma$ is the full tree associated with generation $n$, then $\hFf_n$ will denote the $\sigma$-algebra generated by the $\xi_u$'s with $u$ vertex of generation $k\leq n-1$ and by the $\lambda_v$'s with $v$ vertex of generation $k\leq n$. With these notations, the following result holds

\begin{proposi}
\label{embed11.prop-hausdMartin}
With the assumptions made previously on the construction of the random rooted weighted tree $(\ts,\vs,\bullet,\kappa)$ the following results hold

\noindent (i) For all $s\geq 0$, the family $Y_n(s)=\Hh_n^s/m^nh(s)^n$ is a martingale with respect to the family $\left(\hFf_n\right)_{n\in\NM}$ of $\sigma$-algebras. In particular it converges almost surely to a positive random variable $Y(s)$ such that $\EM(Y(s))=1$.

\noindent (ii) There is $t_m>s_m$ defined as the unique solution of $h(2s)=mh(s)^2$, such that 

(a) if $s<t_m$ then $Y_n(s)$ converges almost surely to a constant, 

(b) if $s=t_m$, then, 

$$\var(Y(t_m)) =
   1-\frac{1}{m}+ \frac{\sigma^2}{m^2}\,,
$$

(c) if $s>t_m$ the random variable $Y(s)$ has not a finite second moment.

\noindent (iii) If $\rho\{1\}<m^{-1}$ and if $s_m$ is the unique solution of $mh(s)=1$, then the sequence of random variables $Y(\Gamma)= \Hh^{s_m}(\Gamma)$ defines a martingale with respect to the family $\left(\hFf_\Gamma\right)_{\Gamma\in\GG}$ of $\sigma$-algebras. In particular it converges almost surely to $1$.
\end{proposi}

\begin{proof} (i) From the definition of $\Hh_n^s$ in eq.~(\ref{embed11.eq-hausmeas}), it follows that

$$\EM\left(\Hh_1^s\right)=
   \EM\left(\sum_{v\in\ch(\bullet)}\lambda^s\right)= 
    h(s) \EM(\xi_\bullet)= mh(s)\,.
$$

\noindent Moreover,

\begin{eqnarray*}
\EM(\Hh_{n+1}^s\,|\,\hFf_n) &=&
  \EM\left(
    \sum_{u\in\vs_{n}}\kappa(u)^s\sum_{v\in\ch(u)} \lambda_v^s\,\Big|\,\hFf_n
     \right)\\
  &=& h(s) \EM\left(
    \sum_{u\in\vs_{n}}\kappa(u)^s\xi_u\,\Big|\,\hFf_n
              \right)
  =m\,h(s) \sum_{u\in\vs_{n}}\kappa(u)^s\\
  &=& mh(s)\Hh_n^s\,.
\end{eqnarray*}

\noindent This calculation shows that $\EM(Y_{n+1}(s)\,|\,\hFf_n)=Y_n(s)$. In particular it is a martingale {\em w.r.t.} the family $\hFf_n$ of $\sigma$-algebras. Since $Y_1(s)=\Hh^s/mh(s)$ it follows that $\EM(Y_1(s))=1$. Therefore $\EM(Y_n(s))=1$ for all $n$'s. The convergence of this family is the main result of the martingale theory \cite{Do53}.

\vspace{.1cm}

\noindent (ii) The calculation of the variance will be done through the second moment of $\Hh_n^s$. By construction

$$\EM\left((\Hh_{n+1}^s)^2\,|\,\hFf_n\right) =
\EM\left(
    \sum_{u,u'\in\vs_{n}}\kappa(u)^s\kappa(u')^s
     \sum_{v\in\ch(u)} \sum_{v'\in\ch(u')} 
      \lambda_v^s\lambda_{v'}^s\,\Big|\,\hFf_n
     \right)\,.
$$

\noindent Let the terms with $u\neq u'$ be considered first. Then, since the $\lambda_v$'s are independent for different $v$'s, it follows that

\begin{eqnarray*}
\EM\left(
    \kappa(u)^s\kappa(u')^s
     \sum_{v\in\ch(u)} \sum_{v'\in\ch(u')} 
      \lambda_v^s\lambda_{v'}^s\,\Big|\,\hFf_n
     \right) &=&
   h(s)^2 \EM\left(
    \kappa(u)^s\kappa(u')^s
     \xi_u\xi_{u'}\,\Big|\,\hFf_n
     \right)\\
&=& h(s)^2 m^2\kappa(u)^s\kappa(u')^s\,.
\end{eqnarray*}

\noindent If now $u=u'$, this gives two terms: the first one are terms for which $v\neq v'$ and the other ones are for $v=v'$. The same type of calculation leads to

\begin{eqnarray*}
\EM\left(
    \kappa(u)^{2s}
     \sum_{v\neq v'\in\ch(u)} 
      \lambda_v^s\lambda_{v'}^s\,\Big|\,\hFf_n
     \right)
&=&
   h(s)^2 (m^2+\sigma^2-m)\,\kappa(u)^{2s}\,,\\
\EM\left(
    \kappa(u)^{2s}
     \sum_{v\in\ch(u)} 
      \lambda_v^{2s}\,\Big|\,\hFf_n
     \right)
&=& h(2s)\, m\,\kappa(u)^{2s}\,.
\end{eqnarray*}

\noindent Grouping these results, leads to

$$\EM\left((\Hh_{n+1}^s)^2\,|\,\hFf_n\right) =
   h(s)^2\,m^2\, (\Hh_n^s)^2+ 
   \{m\left(h(2s)-h(s)^2\right)+\sigma^2\,h(s)^2\}\Hh_n^{2s}\,.
$$

\noindent Averaging on both sides gives

$$\var(\Hh_{n+1}^s) = 
    \left\{m\left(h(2s)-h(s)^2\right)+\sigma^2\,h(s)^2\right\}
     m^{n}h(2s)^n\,.
$$

\noindent It is worth noticing that, thanks to the definition of $h$, the Cauchy-Schwarz inequality gives $h(s)^2 < h(2s)$. This inequality is actually strict because $\rho\{1\}\neq 1$, so that $\lambda$ is not almost surely equal to one. Therefore

$$\var(Y_n(s))= 
  \left(
 \frac{1}{m}\frac{h(2s)}{h(s)^2}
  \right)^n\; 
  \left(
 1-(1-\frac{\sigma^2}{m})\frac{h(s)^2}{h(2s)}
  \right)\,.
$$

\noindent It follows that, if $s<s_m$, then $h(s)m>1$. In addition an elementary calculation shows that the map $g(s) = h(2s)/h(s)^2$ is monotone increasing, that $g(0)=1$ and $\lim_{s\to\infty}g(s)=\rho\{1\}^{-1}$. Therefore, there is a unique $t_m>0$ such that $m=g(t_m)$. Using the definition of $s_m$, it is easy to show that $s_m<t_m$. Hence

(a) if $s<t_m$, $\lim_{n\to\infty} \var(Y_n(s))=0$, implying that $Y_n$ converges almost surely to a constant; this constant can only be the common average, namely $\lim_{n\to\infty}Y_n(s)=1$,

(b) if $s=t_m$, then the variance converges to a finite value

$$s=t_m\;\Rightarrow\;
   \lim_{n\to\infty} \var(Y_n(t_m))=
    \left(
 1+\frac{\sigma^2}{m^2}-\frac{1}{m}
  \right)\,,
$$

(c) if $s>t_m$, then the limiting random variable $Y(s)$ does not have a finite second moment.

\vspace{.1cm}

\noindent (iii) Let now $\Gamma'\subset \Gamma$ be two full finite subtrees of $\ts$. Then there is a decreasing sequence of full finite subtrees such that $\Gamma'\subset \Gamma_j\subset \cdots \Gamma_1\subset \Gamma_0=\Gamma$, and such that $\Gamma_{i+1}$ is obtained from $\Gamma_i$ by the following procedure: each vertex $v\in\Ll\Gamma_i$ which is not in $\Ll\Gamma'$ is removed and replaced by its parent. It is clear that, if $\Gamma_i$ is full, so is $\Gamma_{i+1}$. This leads to 

$$\EM\left(\Hh^s(\Gamma)\,|\, \hFf_{\Gamma_1}\right)=
   \EM\left(
    \sum_{u\in\Ll\Gamma_1;\ch(u)\cap\Ll\Gamma'=\emptyset}
     \kappa(u)^s
      \sum_{v\in\ch(u)} \lambda_v^s\,\Big|\, \hFf_{\Gamma_1}
      \right)+
    \sum_{u\in\Ll\Gamma\cap\Ll\Gamma'}
     \kappa(u)^s
$$

\noindent Thanks to the definition of $\hFf_{\Gamma_1}$ the {\em r.h.s.} becomes 

$$\EM\left(\Hh^s(\Gamma)\,|\, \hFf_{\Gamma_1}\right)=
   mh(s) \sum_{u\in\Ll\Gamma_1;\ch(u)\cap\Ll\Gamma'=\emptyset}
     \kappa(u)^s + \sum_{u\in\Ll\Gamma\cap\Ll\Gamma'}
     \kappa(u)^s\,.
$$

\noindent In particular, if $s=s_m$, namely if $mh(s)=1$, this gives $Y(\Gamma)=\Hh^{s_m}(\Gamma)$ so that

$$\EM\left(Y(\Gamma)\,|\, \hFf_{\Gamma_1}\right)= Y(\Gamma_1)\,.
$$

\noindent Proceeding inductively along the chain of $\Gamma_i$'s, this gives

$$\EM\left(Y(\Gamma)\,|\, \hFf_{\Gamma'}\right)= Y(\Gamma')\,.
$$

\noindent Therefore the family $\left(Y(\Gamma)\right)_{\Gamma\in\GG}$ is also a martingale {\em w.r.t.} the $\hFf_\Gamma$'s. The martingale theorem then shows that it converges almost surely. Since the full tree with boundary $\vs_n$ is a member of this family and since it has been shown that the variance converges to zero (because $s_m<t_m$), the family converges to a constant almost surely.
\end{proof}

\begin{proposi}
\label{embed11.prop-GWHdim}
Under the hypothesis of Proposition~\ref{embed11.prop-hausdMartin}, the Hausdorff dimension of $(\partial \ts,d_\kappa)$ is almost surely equal to $s_m$.
\end{proposi}

\begin{proof} Thanks to Proposition~\ref{embed11.prop-hausdMartin}, $\Hh^{s_m}(\Gamma)$ converges almost surely to $1$. It follows that 

$$\Hh_\delta^{s_m}= 
   \inf_{\Gamma\in\GG_\delta} \Hh^{s_m}(\Gamma)\,,
\hspace{1cm} \Rightarrow\hspace{1cm}
\lim_{\delta\downarrow 0}\Hh_\delta^{s_m}=1\,.
$$

\noindent Consequently, $\Hh_\delta^{s}\to\infty$ for $s<s_m$ and $\Hh_\delta^{s}\to 0$ for $s>s_m$. Hence $\dim_\mathrm{H} (\partial\ts, d_\kappa)=s_m$.
\end{proof}

\begin{proposi}
\label{embed11.prop-GWHmeas}
Under the hypothesis of Proposition~\ref{embed11.prop-hausdMartin}, the Hausdorff measure of $(\partial\ts,d_\kappa)$ at the dimension $s=s_m$ exists almost surely and is a random probability.
\end{proposi}

\begin{proof} In order to prove it, it is sufficient to consider the basis of clopen sets of the form $[u]$ for $u\in\vs$. It boils down to consider 

$$\Hh^{s_m}(\Gamma;u) =
   \sum_{v\in\Ll\Gamma;v\preceq u}
    \kappa(v)^s\,.
$$

\noindent a calculation similar to the one made in the proof of Proposition~\ref{embed11.prop-hausdMartin}, shows that the family of $\{\Hh^{s_m}(\Gamma;u)\,;\, \Gamma\in\Gg\,,\, u\in\vs_\Gamma\}$ is also a martingale satisfying

$$\EM\left(\Hh^{s_m}(\Gamma;u)\,\big|\, \hFf_{\Gamma_0}\right)=
   \kappa(u)^{s_m}\,,
$$

\noindent for all full finite subtree $\Gamma_0$ with $u\in\Ll\Gamma_0$ and $\Gamma\supset \Gamma_0$. Therefore, the Martingale Theorem implies that $\mu([u])=\lim_{\Gamma} \Hh^{s_m}(\Gamma;u)$ exists and that it is a random variable. Since the set of vertices is countable, the set of probability zero on which the convergence does not hold can be chosen independently on $u\in\vs$. By construction

$$\sum_{u\in\Ll\Gamma_0} \Hh^{s_m}(\Gamma;u)=
   \Hh^{s_m}(\Gamma)\,,
$$

\noindent showing that, after taking the limit, $\sum_{u\in\Ll\Gamma_0}\mu([u])=1$. In addition, whenever $\Gamma$ is a full finite subtree, $\{ [v] \ ; \ v \in \Ll \Gamma_0 \}$ forms a partition by clopen sets. Moreover, clopen sets of the form $[v]$ ($v$ a vertex of the tree) generate the $\sigma$-algebra of Borel sets. Hence $\mu$ defines a probability measure on $\partial\ts$. 
\end{proof}



\end{document}